\def\mode{1} 

\if0\mode
\documentclass[journal, twoside, web]{ieeecolor}
\usepackage{generic}
\else
\documentclass[journal,twoside]{IEEEtran}
\fi

\usepackage[english]{babel}
\usepackage[utf8]{inputenc}
\usepackage[T1]{fontenc}
\usepackage{microtype}
\usepackage{setspace}
\usepackage{blindtext}
\usepackage{siunitx}
\usepackage{comment}
\usepackage{xargs}
\usepackage{titlecaps}
\Addlcwords{or with if of and an to via -off through on off for in into the versus vs subject} 

\usepackage[justification=raggedright,singlelinecheck=false,size=footnotesize]{caption}
\usepackage[justification=raggedright,singlelinecheck=false,size=footnotesize]{subcaption}
\usepackage{graphicx}
\usepackage{fancyhdr} 
\usepackage{color}
\usepackage{xcolor}
\usepackage{epsfig}
\graphicspath{{Images/}}

\usepackage{array}
\usepackage{stfloats}
\usepackage{algorithmic}
\usepackage{textcomp}
\let\labelindent\relax
\usepackage{enumitem}
\usepackage[absolute]{textpos}

\usepackage{amsmath,amssymb,amsfonts,amsthm}
\usepackage{mathdots}
\usepackage{bbm,dsfont}
\usepackage{relsize}
\usepackage{nicefrac}
\usepackage{bbm,bm}
\usepackage{mathtools}
\usepackage[mathscr]{eucal}
\usepackage{optidef}

\usepackage{url}
\usepackage[colorlinks,citecolor=blue,linkcolor=blue]{hyperref}
\usepackage[capitalize,nameinlink]{cleveref}
\makeatletter
\newcommand*{\Crefns}[1]{{\@cref@sortfalse\@cref@compressfalse\Cref{#1}}}
\makeatother
\usepackage{orcidlink}
\newcommand\orcidicon[1]{\href{https://orcid.org/#1}{\includegraphics[scale=0.04]{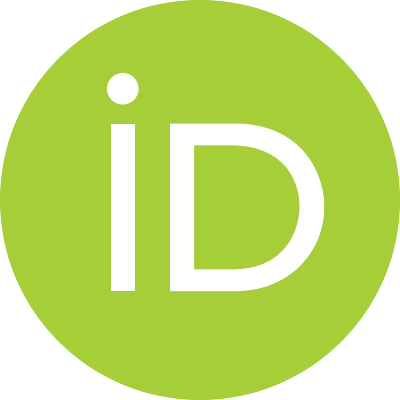}}}

\newcommand{\isExtended}[2]{#1} 

\newcommand{\lr}{\left(}
\newcommand{\rr}{\right)}
\newcommand{\ls}{\left[}
\newcommand{\rs}{\right]}

\newcommand{\Real}[1]{ { {\mathbb R}^{#1} } }

\newcommand{\Complex}[1]{ { {\mathbb C}^{#1} } }

\newcommand{\DeclareAutoPairedDelimiter}[3]{%
	\expandafter\DeclarePairedDelimiter\csname Auto\string#1\endcsname{#2}{#3}%
	\begingroup\edef\x{\endgroup
		\noexpand\DeclareRobustCommand{\noexpand#1}{%
			\expandafter\noexpand\csname Auto\string#1\endcsname*}}%
	\x}

\DeclareAutoPairedDelimiter{\ceil}{\lceil}{\rceil}
\DeclareAutoPairedDelimiter{\floor}{\lfloor}{\rfloor}
\newcommand{\norm}[2][]{\left\lVert#2\right\rVert^2_{#1}}
\newcommand{\normsc}[2][]{\left\lvert#2\right\rvert^2_{#1}}
\newcommand{\e}{\mathrm{e}}
\newcommand{\de}{\mathrm{d}}
\newcommand{\fourier}[1]{\mathcal{F}\ls#1\rs(\lambda)}
\makeatletter
\newcommand{\ifourier}[2][\@empty]{\mathcal{F}^{-1}\ls#2\rs%
    \if#1\@empty \else (#1) \fi}
\makeatother
\DeclareMathOperator*{\erf}{erf}

\newcommand{\spacedom}{\mathcal{X}}
\newcommand{\state}[2]{\psi\!\lr#1,#2\rr}
\newcommand{\statec}[1]{\psi\!\lr#1\rr}
\newcommand{\noise}[2]{v(#1,#2)}
\newcommand{\noisec}[1]{v(#1)}
\newcommandx{\statef}[1]{\hat{\psi}_{\lambda}%
	\ifx#1=={} \else \!\lr#1\rr \fi}
\newcommandx{\noisef}[1]{\hat{v}_{\lambda}%
	\ifx#1=={} \else \!\lr#1\rr \fi}
\newcommand{\ddt}[2]{\dfrac{\partial #1}{\partial t}(#2)}
\newcommand{\dddx}[3]{\dfrac{\partial^2 #1}{\partial x^2}(#2,#3)}
\newcommand{\Aop}{\mathcal{A}}
\newcommand{\Aopf}[1][\lambda]{\hat{A}_#1}
\newcommand{\A}[1]{\ls\Aop#1\rs\!\lr x,t\rr}
\newcommandx{\Aopftot}[2][1={},2=\lambda]{\hat{A}_#2^{#1}}

\newcommand{\delay}{T}
\renewcommand{\u}[2]{u\!\lr#1,#2\rr}
\newcommandx{\uf}[1]{\hat{u}_{\lambda}%
	\ifx#1=={} \else \!\lr#1\rr \fi}
\newcommand{\Kop}[1]{\mathcal{K}_{#1}}
\newcommand{\Kc}[2]{\ls\Kop{}#1\rs\!\lr x,t-#2\rr}
\newcommandx{\Kopk}[2][1={},2={}]{K_{#2}^{#1}}

\newcommandx{\Kopfc}[2][1={},2=\lambda]{\hat{K}_#2%
	\ifx#1=={} \else ^{#1} \fi}

\newcommand{\ub}{\hat{K}_\lambda^\text{u}}
\newcommand{\z}[2]{z(#1,#2)}
\newcommandx{\zf}[1]{\hat{z}_{\lambda}%
	\ifx#1=={} \else \!\lr#1\rr \fi}

\newcommand{\timefreq}{\omega}
\newcommand{\timefreqdom}{\Omega}

\newcommand{\opt}{\text{opt}}

\if0\mode
\newtheoremstyle{colorthm}{\topsep}{\topsep}{\normalfont}{1em}{\color{nblue}\sf\itshape\small\bfseries}{:}{ }{\thmname{#1}\thmnumber{ #2}{\thmnote{ (#3)}}}
\theoremstyle{colorthm}
\fi

\newtheorem{thm}{Theorem}

\newtheorem{prop}{Proposition}
\newtheorem{lemma}{Lemma}

\newtheorem{definition}{Definition}

\newtheorem{ass}{Assumption}
\newtheorem{rem}{Remark}
\newtheorem{ex}{Example}

\Crefname{thm}{Theorem}{Theorems}
\Crefname{cor}{Corollary}{Corollaries}
\Crefname{prop}{Proposition}{Propositions}
\Crefname{conj}{Conjecture}{Conjectures}
\Crefname{ass}{Assumption}{Assumptions}
\Crefname{figure}{Figure}{Figures}
\Crefname{ex}{Example}{Examples}

\addto\captionsenglish{\renewcommand{\figurename}{Fig.}}
\addto\captionsenglish{\renewcommand{\tablename}{TABLE}}

\newcommand{\red}[1]{#1}

\newcommand{\blue}[1]{{\color{blue} #1}}
\newcommand{\juncal}[1]{} 

\newcommand{\revision}[2][]{#1#2}
\newcommand{\linkToPdf}[1]{\href{#1}{\blue{(pdf)}}}
\newcommand{\linkToPpt}[1]{\href{#1}{\blue{(ppt)}}}
\newcommand{\linkToCode}[1]{\href{#1}{\blue{(code)}}}
\newcommand{\linkToWeb}[1]{\href{#1}{\blue{(web)}}}
\newcommand{\linkToVideo}[1]{\href{#1}{\blue{(video)}}}
\newcommand{\linkToMedia}[1]{\href{#1}{\blue{(media)}}}
\newcommand{\award}[1]{\xspace} 
\newcommand{\eg}{\emph{e.g.,}\xspace}
\newcommand{\ie}{\emph{i.e.,}\xspace}
\if0\mode
\addto\extrasenglish{}
\addto\extrasenglish{}
\else
\addto\extrasenglish{}
\addto\extrasenglish{}
\fi

\addto\extrasenglish{}
\addto\extrasenglish{}
\addto\extrasenglish{}

\newcommand{\ptitle}{The role of communication delays in the optimal control of spatially invariant systems}

\title{\titlecap{\ptitle}}
\author{Luca~Ballotta~{\orcidicon{0000-0002-6521-7142}},
	Juncal~Arbelaiz~{\orcidicon{0000-0002-9688-2102}},~\IEEEmembership{Member,~IEEE}, %
	Vijay~Gupta~{\orcidicon{0000-0001-7060-3956}},~\IEEEmembership{Fellow,~IEEE}, %
	Luca~Schenato~{\orcidicon{0000-0003-2544-2553}},~\IEEEmembership{Fellow,~IEEE}, %
	and~Mihailo~R.~Jovanovi\'c~{\orcidicon{0000-0002-4181-2924}},~\IEEEmembership{Fellow,~IEEE}
	\thanks{This work has been partly supported
		by the Italian Ministry of Education, University and Research through
		the PRIN project no. 2017NS9FEY.}%
	\thanks{Luca Ballotta is with the Delft Center for Systems and Control, Delft University of Technology, 2628 CD Delft, Netherlands
		(e-mail: l.ballotta@tudelft.nl).}%
	\thanks{Juncal Arbelaiz is with the Center for Statistics and Machine Learning, Princeton University, Princeton, NJ 08544 USA
		(e-mail: jarbelaiz@schmidtsciencefellows.org).}%
	\thanks{Vijay Gupta is with the Elmore Family School of Electrical and Computer Engineering, Purdue University, West Lafayette, IN 47907 USA
		(e-mail: gupta869@purdue.edu).}%
    \thanks{Luca Schenato is with the Department of Information Engineering, University of Padova, 35131 Padova, Italy
		(e-mail: schenato@dei.unipd.it).}%
	\thanks{Mihailo R. Jovanovi\'c is with the Ming Hsieh Department of Electrical and Computer Engineering, University of Southern California, Los Angeles, CA 90089 USA
		(e-mail: mihailo@usc.edu).}
	}

\begin{document}
	
	\if1\mode
	\begin{textblock}{20}(-2,0.05)
		\footnotesize
		\centering
		\setstretch{1}
		This article has been accepted for publication in the IEEE Transactions on Automatic Control.
		Please cite the paper as:\\
		L. Ballotta, J. Arbelaiz, V. Gupta, L. Schenato, and M. R. Jovanovi\'c,\\
		``\titlecap{\ptitle},''\\
		IEEE Transactions on Automatic Control, 2025.
	\end{textblock}
	\fi
	
	\bstctlcite{MyBSTcontrol}
	
	\maketitle
	\isExtended{}{\numberwithin{equation}{section}}
	

\begin{abstract}
	We study optimal proportional feedback controllers for spatially invariant systems when the controller has access to delayed state measurements received from different spatial locations.
	We analyze how delays affect the spatial locality of the optimal feedback gain leveraging the problem decoupling in the spatial frequency domain.
	For the cases of expensive control and small delay,
	we provide exact expressions of the optimal controllers in the limit for infinite control weight and vanishing delay,
	respectively.
	In the expensive control regime,
	the optimal feedback control law decomposes into a delay-aware filtering of the delayed state and the optimal controller in the delay-free setting.
	Under small delays,
	the optimal controller is a perturbation of the delay-free one which depends linearly on the delay.
 	We illustrate our analytical findings with a reaction-diffusion process over the real line and a multi-agent system coupled through circulant matrices,
 	\revision[showing that delays reduce the effectiveness of optimal feedback control and may require each subsystem within a distributed implementation to communicate with farther-away locations]{}.
	
\end{abstract}

\section{Introduction}\label{sec:intro}

\IEEEPARstart{D}{istributed} and decentralized controller architectures trade performance of the resulting closed-loop system for a reduction in the complexity of the controller implementation~\cite{9117037,ANDERSON2019364}. 
As spatially distributed systems become large-scale, 
this tradeoff turns of utmost importance since the implementation requirements of the all-to-all communication topology inherent to centralized control architectures -- high communication burden, poor scalability, high maintenance costs -- are impractical to sustain. 

The synthesis of optimal distributed controllers requires optimizing the controller in the presence of hard constraints on communications among different subsystems. 
This usually induces non-convex or combinatorial optimization problems that are computationally challenging.
To address such a challenge, 
a popular line of research has been to identify scenarios that lead to tractable optimization problems. 
Examples of such scenarios are spatially invariant systems~\cite{1017552,4623272},
including the case with static proportional feedback and fixed-range communication~\cite{9683472},
sensor selection and statistical modeling~\cite{8876689},
quadratically invariant systems~\cite{XIAO200733,1593901},
and positive systems~\cite{RANTZER201572}.

In this work, we address optimal control of spatially invariant systems,
\revision{namely whose dynamics are translation invariant in the spatial coordinate,}
in the presence of communication delays.
If the open-loop plant is spatially invariant, the optimal delay-free closed-loop system is spatially invariant as well~\cite{1017552}.
\revision{Further, the \linebreak optimal controller is inherently \textit{localized in space} because it implements a static spatial convolution with state measurements, in which the convolution kernel exhibits asymptotic exponential spatial decay~\cite{1017552}.}
This result quantifies the intuition that,
in distributed control,
state measurements from far spatial sites contribute less to the closed-loop performance than measurements from close-by locations.
\revision{Practical implications for the communication design of large-scale systems,
	such as the power grid,
	are significant.
	Long-range links that may be hard to maintain are not actually needed since they add little value to control performance.}

The asymptotic exponential spatial decay bounds have suggested \textit{spatial truncation} of the optimal gain convolution kernel -- \ie \revision{to zero out the feedback gains at locations beyond a fixed distance from the actuator} -- as a heuristic strategy to implement distributed controllers in spatially invariant systems,
\revision{where the cutoff distance relates to, 
	\eg transmission power and bandwidth.}
However,
the optimal kernel is \textit{not} compactly supported in general,
and recovering near-optimal performance through a distributed controller so implemented is not guaranteed under an arbitrary truncation of the optimal kernel. 
In fact, the authors in~\cite{4623272} showed that this truncation can severely degrade the closed-loop performance,
possibly even leading to instability if the truncation is too aggressive. 
A system designer must carefully assess the degree of spatial locality of centralized feedback operators to choose an appropriate ``level of truncation'' and avoid this potential issue~\cite{1017552,4623272,9683472}.

A separate, but related, concern in the design of distributed controllers is that any communication from the subsystems to controllers must occur across non-ideal communication channels. 
We focus here on communication delays that cause a time mismatch between sensing and computation of control input. While the field on how to characterize and mitigate the impact of communication delays on control performance is vast,
see, \eg~\cite{1166540,FRIDMAN2005271,BASIN2007830,8485772,MUNZ20101252,8358743,Michiels2016,5557759}, 
our focus is on characterizing the effect of delays on the spatial locality of optimal feedback gains for spatially invariant systems:
\revision{
	\emph{how should a designer choose feedback gains and communication resources to implement a distributed controller with transmission delays?}}
Relevant works in this realm include~\cite{VOULGARIS2003347,BAMIEH2005575} that studied optimal controllers for cone- and funnel-causal plants where both open-loop dynamics and control actions affect neighboring locations after distance-dependent delays,
and provided structural properties of optimal proportional controllers when the controller has the same cone- or funnel-causal structure of the plant.
The authors in~\cite{FARDAD2011880} extended the previous references to state-space models and tested various optimization algorithms to design distributed controllers.
The recent works~\cite{ballotta2023tcns,Ballotta23lcss-fasterConsensus} report a performance tradeoff between all-to-all and nearest-neighbor controller architectures under architecture-dependent delays,
showing that sparsely connected controllers can be optimal both for mean-square consensus and for fastest convergence in deterministic settings. 
Under similar assumptions for a time-slotted transmission protocol where delays increase with the communication radius,
the authors in~\cite{gupta2010delay} proved that,
while the all-to-all architecture yields the best performance in ring topology,
nearest-neighbor interactions are optimal for plants communicating over regular lattices of dimension greater than one.
Reference~\cite{doi:10.1137/06067897X} targeted the optimal structure for such a controller under specific assumptions.
A related line of work considers nominal,
delay-free feedback and explicitly penalizes the presence of communication links through regularization of a control-theoretic cost function.
Representative examples are~\cite{9117037,8876689},
where sparsity-promoting penalty functions and algorithms are devised to design sparse feedback gains,
the references~\cite{7378905,7347386} that propose communication-aware  \textit{regularization for design},
and the \textit{System Level Synthesis} which imposes local feedback directly on the design of closed-loop maps~\cite{ANDERSON2019364,9939038}.

In this work,
we characterize the impact of communication delays on the spatial locality of optimal feedback controllers for spatially invariant systems.
We study how the suitability of the ``truncation'' approach,
motivated by localization properties of optimal feedback controller kernels for delay-free spatially invariant plants,
is affected by communication delays which are common in spatially distributed settings.
\revision{Differently from early work~\cite{BAMIEH2005575,FARDAD2011880} that focused on plants with finite-speed information flow in both open-loop dynamics and control,}
here we enforce no restrictions on the time evolution of the plant but consider feedback controllers that receive delayed state measurements.
\revision{This is motivated by communication delays present in network systems where different subsystems exchange data over a wireless channel with limited bandwidth.
	While~\cite{BAMIEH2005575,FARDAD2011880} studied the optimal controller in the frequency domain,
	we emphasize how the \emph{spatial} locality of the optimal controller depends on communication delays,}
shedding light on the differences between delay-free and delay-aware optimal controllers for spatially invariant plants; also on how delays affect the suitability of a distributed controller implementation based on spatial truncation of the optimal gain kernel.
\revision{Differently from recent work~\cite{ballotta2023tcns,gupta2010delay},
	we study broader classes of spatially invariant systems not limited to graphs or consensus.}
\revision{The main novelty of this paper,
	which sets it apart from previous works,
	lies in explicitly relating communication delays with spatial locality of optimal feedback gains.
	\revision[This  matters to design both gains and network architecture of a distributed controller,
	and we provide ways to evaluate how a certain amount of delay impacts communication requirements and closed-loop performance compared to delay-free systems]{}.
}

\subsubsection*{Original contribution}

\begin{figure}[t]
	\centering
	\includegraphics[width=.7\linewidth]{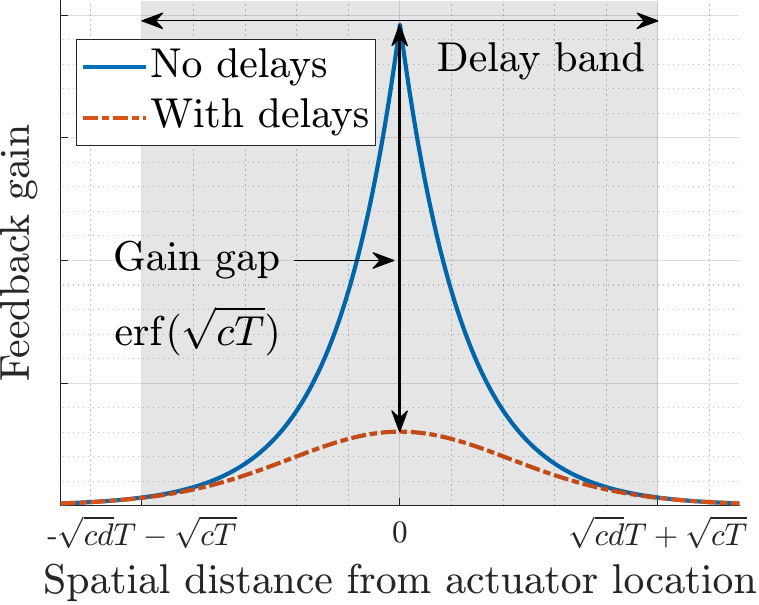}
	\caption{Convolution kernels of optimal feedback controllers for a reaction- diffusion process over the real line. 
		Parameters $c$,
		$d$,
		and $\delay$ are respectively reaction and diffusion coefficients,
		and communication delay.
		\revision[Communication delays decrease and flatten the optimal kernel (dashed red),
		making feedback control less effective and increasing the relative weight of state measurements from distant locations to the actuator.]{}
		In the expensive control regime,
		we analytically quantify the gap between optimal feedback gains at the actuator location without and with delays (``Gain gap'')	and the region where the delay-free and delay-aware gains are substantially different (``Delay band'').
	}
	\label{fig:cover}
\end{figure}

We analytically and numerically investigate the spatial locality of optimal feedback gains for spatially invariant plants subject to temporal delays in the transmission of state measurements used for feedback control.
Since this problem is challenging,
we consider constant delays which nonetheless provide insightful results.
We first analyze the control of \red{scalar} retarded equations,
\revision[concluding that delays make the optimal gains smaller than in the delay-free case;
also, they reduce performance of feedback control]{}.
Building upon this result,
we analytically determine the optimal kernels for spatially invariant plants in two asymptotic cases of practical interest.
In the expensive control regime,
the optimal controller decomposes into a delay-aware filter followed by the optimal controller for delay-free dynamics.
\revision[When truncating the optimal feedback gains to implement a distributed controller,
	the filter may require the retention of measurements from far away that are not used by the delay-free controller and thus increase the communication burden]{}.
Under small delays,
the optimal controller is obtained by perturbing the delay-free kernel with a linear function of the delay,
\revision{which giving indication of the delays the optimal delay-free controller tolerates without overly degrading performance.}
These results are illustrated in two case studies: a reaction-diffusion process and multiple agents coupled through a circulant matrix.

\Cref{fig:cover} depicts the optimal kernels for the former example.
The kernel under delays is flatter than the one without delays;
equivalently, we say that the former has wider spatial spread.
\revision[This implies that the relative value between close-by and far-away measurements for feedback control becomes progressively uniform as a response to increasing delays.
This behavior must be considered in the design of distributed controller implementations via spatial truncation of the gains;
longer delays impose higher communication needs,
such as bandwidth or transmission power,
to exchange measurements with distant subsystems and maintain high control performance.
Our numerical studies corroborate that the optimal gains are smaller due to delays,
which reduces closed-loop performance.]{}

This work considerably expands the paper~\cite{ballotta2023diffusionDelays}
\revision{which studies a reaction-diffusion process in the expensive control regime.}

\subsubsection*{Organization of this article}

\autoref{sec:setup} introduces the main analysis tools, 
including scalar linear retarded systems and the associated optimal control design problem.
In~\autoref{sec:analysis},
we analytically obtain the optimal feedback gains for the asymptotic cases of interest and discuss numerical solutions providing insights into properties of the optimal controllers.
\autoref{sec:centralized} presents our main results for spatially invariant systems.
Optimal controllers for the expensive control regime are characterized in~\autoref{sec:expensive-control-spatial-structure},
and the effect of small delays is described in~\autoref{sec:small-delays-spatially-invariant}.
These two results are further explored through classical dynamical systems.
\autoref{sec:diffusion} considers a reaction-diffusion process with \red{scalar state} evolving on the real line,
which yields the optimal control kernels in~\autoref{fig:cover}.
We analytically approximate the control kernel for the expensive control regime and derive practical design guidelines,
including a discussion on the effects of delays on truncation.
\autoref{sec:circulant} analyzes how small delays affect the closed-loop eigenvalues and the optimal feedback gains for a multi-agent system where agents interact through \red{circulant matrices}.
We draw conclusions and discuss future directions in~\autoref{sec:conclusion}.

\section{Setup}\label{sec:setup}

In this section we lay the basis for the optimal control design problem considered in the rest of the article.
In \autoref{sec:system-model},
we introduce scalar systems with delayed feedback measurements motivated by spatially invariant systems.
In \autoref{sec:problem-formulation},
we formalize the optimal control design problem for such systems.
In \autoref{sec:cost-function},
we provide a closed-form expression of the $\mathcal{H}_2$-norm to be minimized in the design and study the stabilizing feedback gains.
These developments serve as building blocks to design optimal controllers for spatially invariant systems under communication delays,
which is the subject of our main results in \autoref{sec:analysis}. We exemplify in \cref{sec:diffusion,sec:circulant}.

\revision{\subsubsection*{Notation}
	Symbols $\Real{}$ and $\Complex{}$ respectively denote real and complex numbers.
	Square-integrable functions on domain $\spacedom$ are denoted by $L^2(\spacedom)$.
	Operator $\Aop$ has domain $\mathcal{D}(\mathcal{A})$ and maps function $f(x)$ to the image $[\Aop f](x)$.
	Symbols $\mathcal{F}$ and $\mathcal{F}^{-1}$ denote bilateral Fourier transform and inverse Fourier transform,
	respectively.
	Spatial convolution is denoted by $\star$.
	If functions $f$ and $g$ depend on $a\in\Real{}$,
	notation $f \sim g$ means that $\lim_{a\to\ell} \nicefrac{f}{g}=1$ with $\ell$ unambiguously chosen as $\ell\in\{\pm\infty\}$.
}

\subsection{System Model}\label{sec:system-model}

\revision{The scalar linear system we consider is}
\begin{equation}\label{eq:model-spatially-invariant-fourier-control}
	\dfrac{\de\statef{}}{\de t}(t)= \Aopftot\statef{t} +\uf{t} + \noisef{t}
\end{equation}
with state $\statef{t}\in\mathbb{C}$,
control input $\uf{t}\in\mathbb{C}$,
and exogenous disturbance $\noisef{t}\in\mathbb{C}$, where $t \in \mathbb{R}_{\geq 0}$ denotes the instant of time.
Subscript $\lambda$ is a parameter of the system. Its meaning and significance will become clear later on when introducing spatially invariant plants.
Similarly,
the ``hat'' notation,
\eg $\hat{\psi}_\lambda$,
is motivated by the representation of spatially invariant systems in spatial frequency domain,
that we discuss later.
However,
the problem formulation in this section and the analysis in \autoref{sec:analysis} are general and hold for any scalar system.

\revision[Given a well-defined initial value of $\statef{t}$ for $t\in{[-T,0]}$,]{}
we address static feedback from delayed state measurements,
\begin{equation}\label{eq:control}
	\uf{t} = - \Kopfc\statef{t-\delay},
\end{equation}
\revision{where the delay $\delay>0$ is bounded.}
The feedback gain $\Kopfc$ is also parameterized by $\lambda$ and the closed-loop system obtained by applying the control~\eqref{eq:control} to the open-loop plant~\eqref{eq:model-spatially-invariant-fourier-control} is
\begin{equation}\label{eq:model-spatially-invariant-fourier}
	\dfrac{\de\statef{}}{\de t}(t) = \Aopftot\statef{t} - \Kopfc\statef{t-\delay} + \noisef{t}.
\end{equation}
\revision{If the dynamics of a spatially distributed system have added structure in the form of translation invariance in the spatial coordinate (formally defined in \cref{app:preliminaries})},
the spatial Fourier transform maps the system to a family of decoupled subsystems of the form~\eqref{eq:model-spatially-invariant-fourier-control} or~\eqref{eq:model-spatially-invariant-fourier} in spatial frequency domain.
Thus, the abstraction~\eqref{eq:model-spatially-invariant-fourier} is useful to study classes of systems evolving on a spatial domain.
With this interpretation,
$\Aopftot$ and $\Kopfc$ are the Fourier symbols of the open-loop dynamics and controller operators,
respectively; see \cref{app:preliminaries}.
Distribution of such systems across space motivates the time-delay $\delay$ in~\eqref{eq:control} which represents the communication delay incurred when collecting measurements from spatially distant subsystems.
We point the reader to~\cite{1017552} for details.
To ground the discussion,
we introduce two common instances of spatially invariant systems.

\begin{ex}[Reaction-diffusion process]\label{ex:reaction-diffusion}
	The dynamics of a linear reaction-diffusion process are
	\begin{subequations}\label{eq:reaction_diff_joint}
		\begin{equation}\label{eq:model-diffusion-ex}
			\ddt{\psi}{x,t} = d\dddx{\psi}{x}{t} - c\state{x}{t} + \u{x}{t} + \noise{x}{t}.
		\end{equation}
		The state $\state{\cdot}{t}:\Real{}\to\Real{}$ represents a time-varying quantity (\eg heat) distributed along an infinite one-dimensional spatial domain.
		Constants $c>0$ and $d>0$ are reaction and diffusion coefficients,
		respectively.
		We consider feedback controllers computed 
		by spatially convolving  delayed state measurements with the {convolutional kernel} (or {control kernel}) $\Kopk:\Real{}\to\Real{}$:
		\begin{equation}\label{eq:controller-convolution-ex}
			\u{x}{t} = -\dfrac{1}{\sqrt{2\pi}}\int_{\Real{}} \Kopk(x - y)\state{y}{t-\delay}\de y.
		\end{equation}
	\end{subequations}
	\revision{Here,
		the spatial locality of the controller is given by the shape of its kernel.
		If $\Kopk$ is nearly flat,
		namely feedback gains have similar magnitude,
		measurements from far away matter for feedback,
		while they do not if $\Kopk$ features a sharp spatial decay.}
	
	In spatial frequency domain,
	provided that the state $\state{\cdot}{t}$ and the control kernel $K$ have well-defined spatial Fourier transforms~\cite{9683472},
	systems~\eqref{eq:model-diffusion-ex} and~\eqref{eq:reaction_diff_joint} are  mapped to a family of systems of the kind~\eqref{eq:model-spatially-invariant-fourier-control} and~\eqref{eq:model-spatially-invariant-fourier} respectively, parameterized by the spatial frequency $\lambda\in\Real{}$,
	where $\Aopftot=-d\lambda^2-c$ and $\Kopfc$ is the spatial Fourier transform of $\Kopk$.
	We study optimal control of reaction-diffusion processes in~\autoref{sec:diffusion}.
\end{ex}

\begin{ex}[Multi-agent system]\label{ex:circulant}
	Consider the closed-loop dynamics on a finite number of integer-indexed spatial locations
	\begin{equation}\label{eq:circulant-dynamics}
		\dfrac{\de\psi}{\de t}(t) = \Aop\statec{t} - \Kop{}\statec{t-\delay} + \noisec{t}
	\end{equation}
	where $\statec{t} \in \Real{N}$ is a real-valued vector and
	$\Aop\in\Real{N\times N}$ and $\Kop{}\in\Real{N\times N}$ are circulant matrices.
	This model can represent a multi-agent system with $N$ agents where each agent $i\in\{0,\dots,N-1\}$ has a scalar state $\state{i}{t}$ and	$\Kop{}$ is the feedback gain matrix.\footnote{Multiplying the state $\statec{t}$ by a circulant matrix corresponds to a spatial convolution between $\statec{t}$ and the first row of the matrix.}
	Inter-agent couplings in the dynamics are modeled by suitably structuring the open-loop matrix $\Aop$,
	\eg nearest neighbor interactions correspond to a tridiagonal matrix with two additional nonzero elements in the top-right and bottom-left positions.
	The representation of~\eqref{eq:circulant-dynamics} in the spatial frequency domain is given by~\eqref{eq:model-spatially-invariant-fourier} for $\lambda\in\{0,\dots,N-1\}$ where $\Aopftot$ and $\Kopfc$ are the (possibly scaled) eigenvalues of $\Aop$ and $\Kop{}$,
	respectively.
	We analyze this example in~\autoref{sec:circulant}.
\end{ex}

\subsection{Problem Formulation}\label{sec:problem-formulation}

\revision{Given weight $r>0$ that determines how much control effort is penalized by the design,
	consider the performance output}
\begin{equation}\label{eq:performance-spatial-transform}
	\zf{t} \doteq \begin{bmatrix}
		\statef{t}\\
		\sqrt{r}\uf{t}
	\end{bmatrix} = \begin{bmatrix}
		\statef{t}\\
		-\sqrt{r}\,\Kopfc\statef{t-\delay}
	\end{bmatrix}.
\end{equation}
We address $\mathcal{H}_2$-optimal control of the closed-loop system~\eqref{eq:model-spatially-invariant-fourier} with output $\zf{}$.
The transfer function from $\noisef{}$ to $\statef{}$ is
\begin{equation}\label{eq:tf}
	\hat{H}_\lambda(s) \doteq \dfrac{1}{s - \Aopftot + \Kopfc\e^{-\delay s}},
\end{equation}
where $s\in\mathbb{C}$ is the Laplace variable,
and the $\mathcal{H}_2$-norm from $\noisef{}$ to $\zf{}$,
denoted by $J_\lambda$,
is
\begin{equation}\label{eq:h2-norm-decoupled}
	J_\lambda(\Kopfc) = \lr 1+r\normsc{\Kopfc}\rr \int_{\timefreqdom} \normsc{\hat{H}_\lambda(j\timefreq)} \de \timefreq,
\end{equation}
where $\timefreqdom$ is the time frequency domain.
We are interested in the following optimal control design problem:
\begin{mini!}
	{\Kopfc \in \mathbb{C}}
	{J_\lambda(\Kopfc)}
	{\label{eq:prob-optimal-control-decoupled}}
	{}
	\addConstraint{\eqref{eq:model-spatially-invariant-fourier} \mbox{ is stable.}\label{eq:prob-optimal-control-decoupled-constr}}
\end{mini!}
\revision{Stable delay systems converge exponentially fast~\cite{datko1978procedure} and hence the $\mathcal{H}_2$-norm $J_\lambda$ is well defined~\cite{5557759}.
	Problem~\eqref{eq:prob-optimal-control-decoupled} heuristically reduces control energy and state deviation due to disturbances.
	If $\noisef{}$ is white noise,
	then $J_\lambda$ coincides with the expected energy of the performance output $\zf{}$.}
We denote the solution of~\eqref{eq:prob-optimal-control-decoupled} by $\Kopfc[\opt][\delay](\lambda)$ to highlight that it depends on the delay $\delay$ in the measurement transmission. This problem formulation allows us to study optimal controllers for classes of spatially invariant systems such as those in \cref{ex:reaction-diffusion,ex:circulant}.
Similarly to the delay-free setting in~\cite{1017552},
solving~\eqref{eq:prob-optimal-control-decoupled} for all spatial frequencies $\lambda$ 
returns the Fourier symbol of the optimal controller for a spatially invariant system.
The optimal controller in physical domain is retrieved from the inverse spatial Fourier transform of 
the solution to~\eqref{eq:prob-optimal-control-decoupled}.
This corresponds to the kernel $\Kopk[\opt][\delay]$ in \cref{ex:reaction-diffusion} and to the feedback gain matrix $\Kop{\delay}^\opt$ in \cref{ex:circulant}.

\subsection{Cost Function and Stabilizing Feedback Gains}\label{sec:cost-function}

\begin{figure}
    \centering
    \includegraphics[height=.55\linewidth]{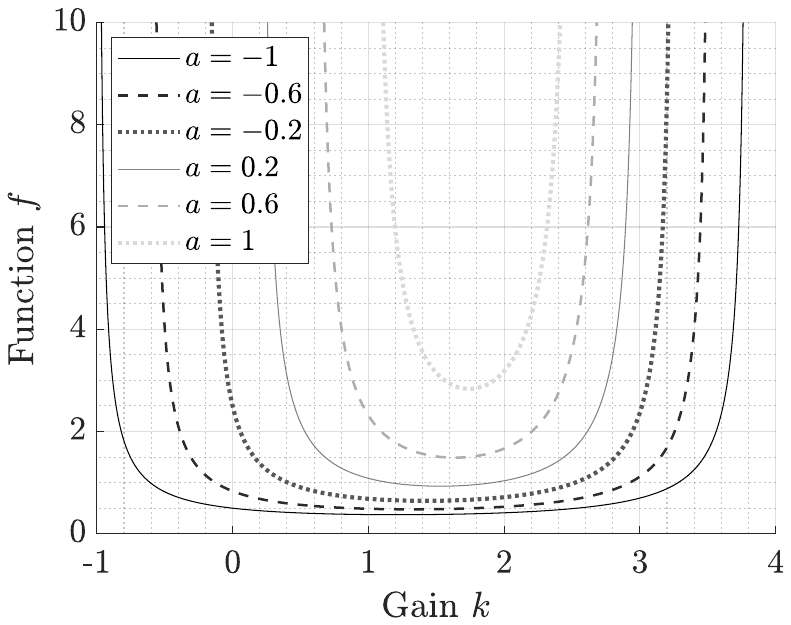}
    \caption{Graphic of $f$ in~\eqref{eq:time-integral-function} with $\delay=0.5$ and $a\in\{\pm0.2,\pm0.6,\pm1\}$.}
    \label{fig:variance}
\end{figure}

In the general case when either $\Aopf$ or $\Kopfc$ in~\eqref{eq:model-spatially-invariant-fourier} is complex,
explicitly computing the cost $J_\lambda$ in~\eqref{eq:prob-optimal-control-decoupled} analytically is challenging.
However, several systems of interest -- such as those provided in \cref{ex:reaction-diffusion,ex:circulant} -- have real Fourier symbols. Thus, we make the following assumption for the remainder of this paper, which provides analytical tractability.

\begin{ass}[Real Fourier symbols]\label{ass:A-even}
	The open-loop coefficient $\Aopftot$ and the feedback gain $\Kopfc$ in~\eqref{eq:model-spatially-invariant-fourier} are real.
\end{ass}
\cref{ass:A-even} implies that the real and imaginary parts of $\statef{t}$ in~\eqref{eq:model-spatially-invariant-fourier} evolve independently.
Hence,
without loss of generality,
in the following we assume that $\statef{t}\in\Real{}$.
\revision{Given \cref{ass:A-even},
	we compute the cost $J_\lambda$ as described next.}

The authors in~\cite{KuchlerLangevinEqs} studied stationary solutions for a real-valued stochastic process $x$ that obeys linear retarded dynamics stochastically forced  by a standard Wiener process $w$,
\begin{equation}\label{eq:sdde}
	\de x(t) = (ax(t) - kx(t-\delay))\de t + \de w(t),
\end{equation}
\revision{where $a,k\in\Real{}$ which requires \cref{ass:A-even}.}
If a stationary solution to~\eqref{eq:sdde} exists,
it is Gaussian with steady-state variance equal to the energy of the so-called fundamental solution associated with~\eqref{eq:sdde}.
The fundamental solution is the signal $x_0$ that obeys the deterministic dynamics obtained by removing $w$ from~\eqref{eq:sdde} and imposing the initial condition $x_0(0) = 1$,
$x_0(t) = 0 \; \forall t \in [-\delay,0)$.
Applying Laplace transform and Parseval's theorem,
the squared $L_2$-norm of $x_0$ amounts to
\begin{equation}\label{eq:time-integral-sdde}
	\int_{0}^{+\infty} x_0^2(t)\de t = \int_{-\infty}^{+\infty} \dfrac{\de\timefreq}{\normsc{j\timefreq - a + k\e^{-j\delay\timefreq}}}.
\end{equation}
The right-hand side of~\eqref{eq:time-integral-sdde} is the integral of $\lvert\hat{H}_\lambda(j\timefreq)\rvert^2$ over frequency domain
after replacing $a $ with $ \hat{A}_\lambda$ and $k$ with $\hat{K}_\lambda$.
Hence,
we use the result in~\cite{KuchlerLangevinEqs} that computes~\eqref{eq:time-integral-sdde} to evaluate the integral in~\eqref{eq:h2-norm-decoupled}.
We denote the integral~\eqref{eq:time-integral-sdde} by $f(k)$,
where we make the dependence on the feedback gain $k$ explicit.
Its closed-form expression derived in~\cite{KuchlerLangevinEqs} \revision{for $aT<1$} is
{\renewcommandx{\Kopfc}[2][1={},2={}]{k_{#2}^{#1}}
	\renewcommand{\Aopf}{a}
	\renewcommand{\Aopftot}[1][]{a^{#1}}
	\renewcommand{\ub}{k_\text{u}}
	\begin{equation}\label{eq:time-integral-function}
		f(\Kopfc) = 
		\begin{cases}
			\dfrac{-\Kopfc\sinh(\ell\delay) - \ell}{2\ell\lr \Aopftot-\Kopfc\cosh(\ell\delay)\rr},	& \lvert \Kopfc\rvert < -\Aopftot\\
			\dfrac{\delay}{4} + \dfrac{1}{4|\Aopftot|},												& \Kopfc = |\Aopftot|, \Aopftot\neq0\\
			\dfrac{-\Kopfc\sin(\ell\delay) - \ell}{2\ell\lr \Aopftot-\Kopfc\cos(\ell\delay)\rr},	& |\Aopftot| < \Kopfc < \ub
		\end{cases}
	\end{equation}
	where $\ell \doteq \sqrt{|\Kopfc[2] - \Aopftot[2]|}$.
	The upper bound on stabilizing gains,
	$\ub$,
	is the unique solution of the following equation for $\Kopfc > |\Aopf|$:
	\begin{equation}\label{eq:upper-bound}
		\delay\sqrt{\Kopfc[2] - \Aopftot[2]} = \arccos\lr\Aopftot\Kopfc[-1]\rr.
	\end{equation}
	The profile of $f$ is shown in~\autoref{fig:variance} for several values of $a$.
	For given $a$ and $\delay$,
	the stabilizing gains belong to the interval $(a,\ub)$,
	which thus determines the stability condition~\eqref{eq:prob-optimal-control-decoupled-constr}.
	\revision{All three cases in~\eqref{eq:time-integral-function} are well defined if $\Aopf < 0$,
		while only the third case is if $\Aopf \ge0$ and $\Aopf\delay<1$.
		The $\mathcal{H}_2$-norm is undefined if $\Aopf\delay\ge1$ because system~\eqref{eq:sdde} or its deterministic counterpart cannot be stabilized under this condition~\cite{KuchlerLangevinEqs,5557759}.}
	The next results characterize the upper bound $k^u$ as $a$ and $\delay$ vary.
	
	\begin{prop}[Stabilizing gains vs. dynamics]\label{prop:upper-bound}
		The upper bound $\ub$ is decreasing with $\Aopf$ and convex in $\Aopf$.
		For $\delay>0$,
		it holds	
		\begin{equation}\label{eq:upper-bound-limit}
			\lim_{\Aopf\rightarrow-\infty}\dfrac{\ub}{|\Aopf|} = 1, \qquad \lim_{\Aopf\rightarrow\frac{1}{\delay}}\ub = \dfrac{1}{\delay}.
		\end{equation}
		\begin{proof}
			See~\cref{app:upper-bound}.
		\end{proof}
	\end{prop}
	\begin{prop}[Stabilizing gains vs. delay]\label{prop:upper-bound-vs-delay}
		The upper bound $\ub$ is decreasing with $\delay$ and convex in $\delay$.
		For $\Aopf>0$,
		it holds
		\begin{subequations}\label{eq:upper-bound-limit-del}
			\begin{equation}\label{eq:upper-bound-limit-del-1}
				\lim_{\delay\rightarrow\frac{1}{\Aopf}}\ub = \Aopf, \qquad 
				\lim_{\delay\rightarrow0^+}\ub = +\infty
			\end{equation}
			and,
			for $\Aopf<0$,
			it holds
			\begin{equation}\label{eq:upper-bound-limit-del-2}
				\lim_{\delay\rightarrow+\infty}\ub = |\Aopf|, \qquad 
				\lim_{\delay\rightarrow0^+}\ub = +\infty
			\end{equation}
		\end{subequations}
		\begin{proof}
			Analogous to the proof of~\cref{prop:upper-bound}.
		\end{proof}
	\end{prop}
	
	\Cref{prop:upper-bound} implies that stabilizing feedback controllers for~\eqref{eq:sdde} are conservative,
	\ie feature ``small'' feedback control gains,
	if the timescale of the open-loop dynamics is comparable to the communication delay in magnitude,
	as suggested by the second limit in~\eqref{eq:upper-bound-limit}.
	Stable but slow open-loop dynamics and unstable plants are stabilized only with small feedback gains,
	\revision[which reduces the effectiveness of control]{}.
	Moreover,
	while all stable open-loop systems have a nonempty set of stabilizing gains,
	unstable plants can be stabilized only under the condition $\Aopf\delay < 1$,
	which means that the delay cannot be larger than the time constant of the system.
	Intuitively,
	the control actions are stabilizing only if the feedback measurements (generated at time $t-\delay$) are sufficiently close in time with the current state of the system at time $t$.
	The insight of \cref{prop:upper-bound-vs-delay} is similar.
	Limits~\eqref{eq:upper-bound-limit-del} suggest that stabilizing gains get arbitrarily large as the delay vanishes,
	\revision[gradually recovering the delay-free case]{},
	whereas for large delays the control gains cannot be greater than the inverse of the time constant of the open-loop system.
}

\section{Optimal Feedback Control Gains}\label{sec:analysis}

Equipped with the exact expression of the $\mathcal{H}_2$-norm, $J_\lambda$,
the optimal control design problem~\eqref{eq:prob-optimal-control-decoupled} under \cref{ass:A-even} is
\begin{mini!}
	{\Kopfc \in \mathbb{R}}
	{J_\lambda(\Kopfc) = \lr 1+r\Kopfc[2]\rr f_\lambda(\Kopfc)\label{eq:prob-optimal-control-decoupled-explicit-obj}}
	{\label{eq:prob-optimal-control-decoupled-explicit}}
	{}
	\addConstraint{\Aopftot}{< \Kopfc<\ub\protect\label{eq:prob-optimal-control-stability-region},}
\end{mini!}
where $f_\lambda$ and $\ub$ respectively denote function~\eqref{eq:time-integral-function} and the upper bound on stabilizing feedback gains $k^\text{u}$ with $\Aopf$ in place of $a$.
The constraint~\eqref{eq:prob-optimal-control-stability-region} corresponds to stable closed-loop systems,
for which the $\mathcal{H}_2$-norm is well defined.

Analytically obtaining the solution of~\eqref{eq:prob-optimal-control-decoupled-explicit} is still challenging due to the involved expression of $f_\lambda$;
in general, the solution needs to be computed numerically.
Moreover,
while numerical tests suggest that $J_\lambda$ is strictly convex,
this is challenging to prove.
We next delve into obtaining analytical solutions of~\eqref{eq:prob-optimal-control-decoupled-explicit} for asymptotic cases of interest, which provide insight into the optimal feedback gains.
In \autoref{sec:tails},
we address stable open-loops dynamics with $\Aopftot\rightarrow-\infty$.
In \autoref{sec:expensive-control},
we examine the expensive control regime $r\rightarrow+\infty$ for stable open-loop dynamics.
In \autoref{sec:small-delays},
we consider vanishing delays $\delay\rightarrow0^+$.
Leveraging the insights from these asymptotic cases,
we present and discuss numerical solutions in \autoref{sec:control-dde-numerical}.

\subsection{\titlecap{Fast stable open-loop plants}}\label{sec:tails}

The next result characterizes the limit behavior of the minimizer of $J_\lambda$ for stable systems	with short timescale.

\begin{prop}[Optimal control gain for fast open-loop dynamics]\label{prop:optimal-gain-tails}
	Under \cref{ass:A-even} and for $\delay>0,r>0$,
	it holds
	\begin{equation}\label{eq:kstar-tails-approx}
		\lim_{\Aopf\rightarrow-\infty}\dfrac{1}{\Kopfc[\opt][\delay](\lambda)}\dfrac{\e^{\delay\Aopftot}}{2r|\Aopftot|} = 1.
	\end{equation}
	\begin{proof}
		See~\cref{app:optimal-gain-small-a}.
	\end{proof}
\end{prop}
\cref{prop:optimal-gain-tails} reveals that the optimal feedback gain features an \textit{asymptotic exponential decrease} with rate $\delay\Aopf$.
In contrast,
the standard delay-free ($\delay=0$) optimal feedback gain
\begin{equation}\label{eq:kstar-no-delay}
	\Kopfc[\opt][0](\lambda) = \Aopftot + \sqrt{\Aopftot[2]{} + \dfrac{1}{r}}
\end{equation}
has,
for $\Aopftot\to-\infty$,
the asymptotic expression
\begin{equation}\label{eq:kstar-no-delay-tails-approx}
	\Kopfc[\opt][0](\lambda) = \dfrac{1}{2r|\Aopftot|} + o\lr\dfrac{1}{|\Aopftot|}\rr,
\end{equation}
which shows that the optimal control feedback gain in the delay-free case is \textit{asymptotically inversely proportional} to $|\Aopftot|$ as this grows unbounded.

We define the characteristic timescale of the open-loop dynamics by $t_* \doteq \nicefrac{1}{|\Aopf|}$. 
\textit{Dimensional analysis} shows that the dimensionless parameter $\alpha_* \doteq \nicefrac{\delay}{t_*}$ -- which captures the ratio of the timescale of measurement delays and the timescale of the open-loop dynamics -- provides physical interpretation of the value of delayed measurements for optimal feedback control~\cite{Arbelaiz2024}.
The exponential $e^{-\alpha_*}$ scales the magnitude of the (asymptotically) optimal feedback gain $\Kopfc[\opt][\delay](\lambda)$,
\revision[reducing the aggressiveness of optimal feedback control compared to the delay-free case with gain {$\Kopfc[\opt][0](\lambda)$}]{}.
Hence,
we enlighten the effects of time delays on the optimal gain by considering the following two regimes.
\begin{description}
	\item[\boldmath Regime $\alpha_* \ll 1$ (\ie $t_* \gg \delay$):] In this case, 
	the time-delay $\delay$ in the feedback control is small compared to the characteristic timescale $t_*$ of the stable open-loop dynamics.
	Thus, 
	the feedback control signal	is useful to optimally stabilize the system despite delayed measurements,
	and the optimal delay-aware gain $\Kopfc[\opt][\delay](\lambda)$ resembles $\Kopfc[\opt][0](\lambda)$.
	\item[\boldmath Regime $\alpha_* \gg 1$ (\ie $t_* \ll \delay$):]  
	In this regime, the open-loop stable dynamics have a timescale that is much shorter than the communication delay,
	so feeding the controller with delayed measurements is useless for optimal stabilization.
	Consequently, 
	$\Kopfc[\opt][\delay](\lambda) \to 0$ to avoid  unnecessary control cost in the performance objective $J_\lambda$.
\end{description}

\subsection{Expensive Control Regime}\label{sec:expensive-control}

We now analyze the setting where the control effort penalization is much higher than penalization in the state deviation.
We consider open-loop stable systems for which it is possible to analytically compute the optimal feedback gain.
\begin{prop}[Optimal control in expensive regime]\label{thm:optimal-gain-expensive-regime}
	Under \cref{ass:A-even} and for fixed $\Aopftot<0$,
	it holds
	\begin{equation}\label{eq:kstar-expensive-approx-lim}
		\lim_{r\rightarrow+\infty}\dfrac{1}{\Kopfc[\opt][\delay](\lambda)}\dfrac{\e^{\delay\Aopftot}}{2r|\Aopftot|} = 1.
	\end{equation}
	\begin{proof}
		See~\cref{app:optimal-gain-expensive-regime}.
	\end{proof}
\end{prop}

\cref{thm:optimal-gain-expensive-regime} implies that,
for large enough $r$ and with stable open-loop dynamics,
the optimal feedback gain $\Kopfc[\opt][\delay](\lambda)$ is approximated by $\Kopfc[\text{ex}][\delay](\lambda)$,
defined next.
We refer to $\Kopfc[\text{ex}][\delay](\lambda)$ as
the \textit{optimal gain in the expensive control regime}:
\begin{equation}\label{eq:kstar-expensive-approx}
	\Kopfc[\opt][\delay](\lambda) \sim \Kopfc[\text{ex}][\delay](\lambda) \doteq \dfrac{\e^{\delay\Aopftot}}{2r|\Aopftot|}.
\end{equation}
The notation $\sim$ means that $\Kopfc[\opt][\delay](\lambda)$ and $\Kopfc[\text{ex}][\delay](\lambda)$ are asymptotically equal at the limit for $r\rightarrow\infty$,
according to~\eqref{eq:kstar-expensive-approx-lim}.
Notably,
the same asymptotic expression holds in the limit for short timescale $\Aopftot\rightarrow-\infty$ according to \cref{prop:optimal-gain-tails}.
We compare the approximation~\eqref{eq:kstar-expensive-approx} to the delay-free case by a Taylor expansion of~\eqref{eq:kstar-no-delay},
which yields the same asymptotic expression in~\eqref{eq:kstar-tails-approx} (with $\nicefrac{1}{r}$ inside the little-o).
Analogously to what discussed for the short timescale regime in the previous section,
the expression of $\Kopfc[\text{ex}][\delay](\lambda)$ in~\eqref{eq:kstar-expensive-approx} agrees with $\Kopfc[\text{ex}][0](\lambda)$ in~\eqref{eq:kstar-tails-approx}
except for the factor $\e^{\delay\Aopftot}$ that implements a correction for delays.
Given the assumption $\Aopftot < 0$,
this correction factor is smaller than one and strictly decreases with $\delay$.
This is consistent with the intuition that delays force smaller controller gains and thus more conservative control actions.

\subsection{Small Delays}\label{sec:small-delays}

We now focus on the setting where the delay $\delay$ is negligible compared to the characteristic timescale of the open-loop dynamics~\eqref{eq:model-spatially-invariant-fourier},
that is, the regime $\alpha_*\ll1$ discussed in \autoref{sec:tails}.
In this case,
it runs out that \red{the optimal controller is a perturbed version of the delay-free controller,
where the perturbation is proportional to the delay.}

\begin{prop}[Optimal control with small delays]\label{prop:optimal-gain-small-delay}
	Under \cref{ass:A-even},
        it holds:
	\begin{equation}\label{eq:kstar-small-T}
		\lim_{\delay\rightarrow0^+} \dfrac{\Kopfc[\opt][\delay](\lambda)}{\Kopfc[\opt][0](\lambda) - \lr\Aopftot\Kopfc[\opt][0](\lambda) + \dfrac{1}{r}\rr\delay} = 1.
	\end{equation}
	\begin{proof}
		See~\cref{app:optimal-gain-small-delay}.
	\end{proof}
\end{prop}
\cref{prop:optimal-gain-small-delay} implies that the optimal feedback gain can be approximated as follows for small enough delay $\delay$:
\begin{equation}\label{eq:kstar-small-T-approx}
	\Kopfc[\opt][\delay](\lambda) \sim \Kopfc[\text{sd}][\delay](\lambda) \doteq \underbrace{\Kopfc[\opt][0](\lambda)}_{\substack{\text{optimal} \\ \text{delay-free gain}}} - \underbrace{\lr\Aopftot\Kopfc[\opt][0](\lambda) + \dfrac{1}{r}\rr\delay}_{\text{correction term}}.
\end{equation}
The feedback gain $\Kopfc[\text{sd}][\delay](\lambda)$ of the optimal controller under small delays is approximated by 
the difference between the optimal delay-free gain and a corrective term which is linear with the delay $\delay$.
Plugging in~\eqref{eq:kstar-no-delay} shows that
for every bounded $\Aopf$
\begin{equation}\label{eq:kstar-small-T-general-correction-term}
	\Aopftot\Kopfc[\opt][0](\lambda) + r^{-1} > 0.
\end{equation}
This makes the correction term in~\eqref{eq:kstar-small-T-approx} negative and subsequently,
for $\delay$ small enough,
smaller gain $\Kopfc[\opt][\delay](\lambda)\le\Kopfc[\opt][0](\lambda)$ than the delay-free optimal controller,
similarly to the two regimes considered in \cref{sec:tails,sec:expensive-control}.

\subsection{Numerical Solutions}\label{sec:control-dde-numerical}

\begin{figure}
	\centering
	\begin{subfigure}{\linewidth}
		\centering
		\includegraphics[width=.7\linewidth]{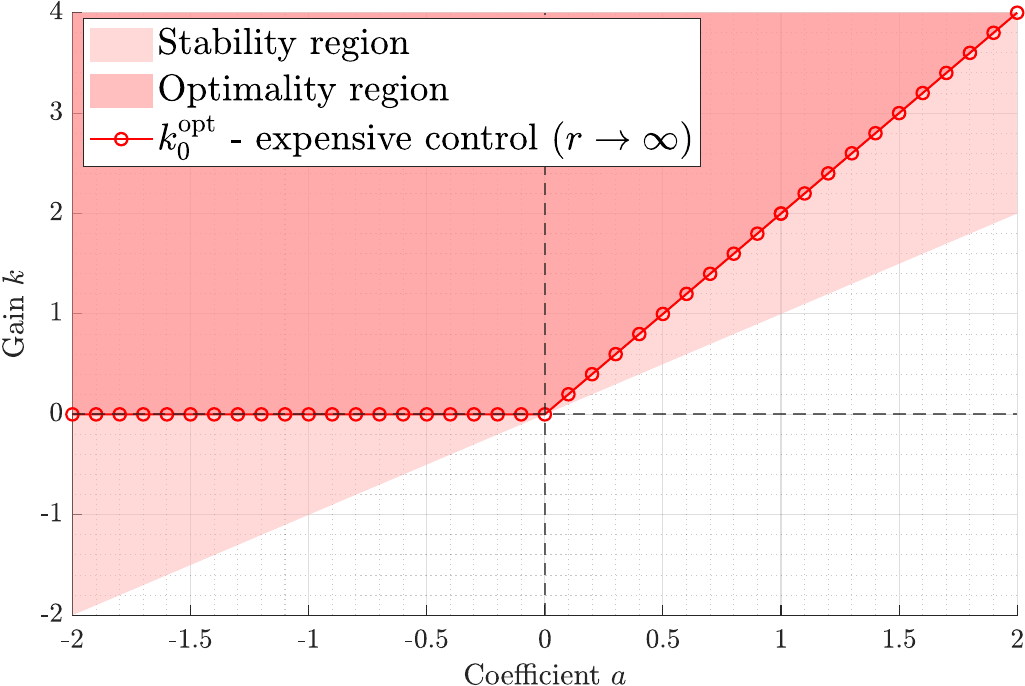}
		\caption{Delay-free case, \ie $\delay=0$ in~\eqref{eq:control}.}
		\label{fig:stab-region-no-del}
	\end{subfigure}
	\hfill
	\vspace{1mm}
	\begin{subfigure}{\linewidth}
		\centering
		\includegraphics[width=.7\linewidth]{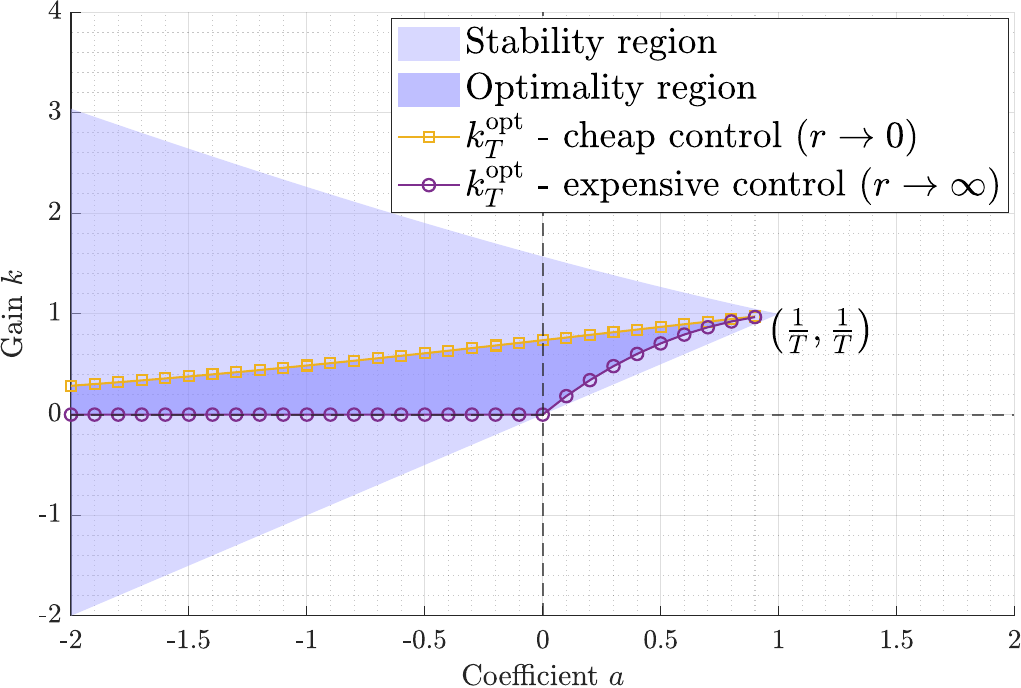}
		\caption{Delay $\delay=1$ in~\eqref{eq:control}.}
		\label{fig:stab-region-del}
	\end{subfigure}
	\caption{Stability (light colored) and optimality (dark colored) regions for optimal control design problem~\eqref{eq:prob-optimal-control-decoupled-explicit}.}
	\label{fig:stab-region}
\end{figure}

We numerically solve problem~\eqref{eq:prob-optimal-control-decoupled-explicit} to complement the analytic asymptotic results derived so far.
We compare the delayed feedback setting under consideration with the standard delay-free one that features $\delay=0$ in~\eqref{eq:control},
for which the $\mathcal{H}_2$-norm coincides with the limit of $f_\lambda$ as $\delay\rightarrow0^+$.

\revision{The ``stability region'' in~\autoref{fig:stab-region-del} depicts the interval~\eqref{eq:prob-optimal-control-stability-region} as $\Aopf$ varies.
	It represents all and only gains that stabilize the closed-loop dynamics.
	\Cref{fig:stab-region-no-del} shows the analogous region for the delay-free system}.
Contrarily to the delay-free case,
the stabilizing gains are upper-bounded in the presence of delays and exist only if $\Aopf\delay < 1$.
The upper bound $\ub$ decreases with $\Aopf$ and tends to $\nicefrac{1}{\delay}$ as $\Aopf\rightarrow\nicefrac{1}{\delay}$,
according to \cref{prop:upper-bound}.

\revision{The ``optimality region'' in \autoref{fig:stab-region-del} delimits the ``meaningful'' portion of the stability region since it contains the solutions of~\eqref{eq:prob-optimal-control-decoupled} for all choices of the control cost parameter $r$.
	The rest of the stability region is never optimal for the $\mathcal{H}_2$-norm~\eqref{eq:h2-norm-decoupled}.
	This relates to the inverse optimal control problem which searches for the cost parameter $r$ that makes a specific choice of the gain $\Kopfc$ optimal~\cite{kalmanInverseOptimality}.}
The boundaries of the (open) optimality region are the optimal control gains in the cheap and expensive control regimes,
corresponding to the limit cases $r\rightarrow0$ and $r\rightarrow+\infty$,
respectively.
We numerically compute the boundaries in~\autoref{fig:stab-region-del} by setting $J(\Kopfc) = f(\Kopfc)$ and $J(\Kopfc) = \Kopfc[2] f(\Kopfc)$ in~\eqref{eq:prob-optimal-control-decoupled-explicit},
respectively.
The delay fundamentally \revision[reduces]{} the optimality region compared to the delay-free setting,
\revision[limiting the effectiveness of feedback control]{}.
With delays,
the optimal feedback gains are upper bounded for every $r>0$ and decrease to zero as the open-loop coefficient $\Aopf$ decreases to negative infinity,
according to~\cref{prop:optimal-gain-tails}.
In contrast, 
the optimal feedback gains in the delay-free setting grow unbounded as $r$ decreases for every $\Aopf$,
see \autoref{fig:stab-region-no-del},
\revision[which makes aggressive control optimal even with stable open-loop dynamics]{}.

\begin{figure}
	\centering
	\includegraphics[width=.7\linewidth]{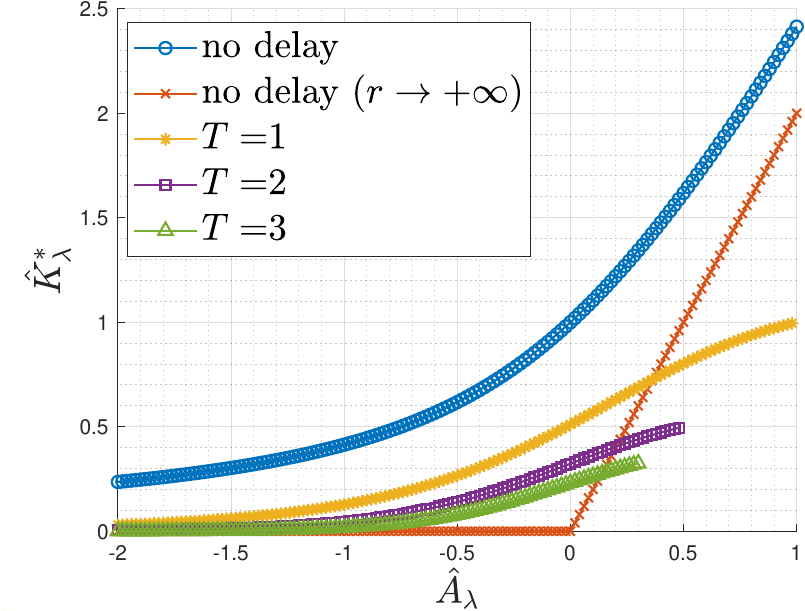}\caption{Optimal gain $\Kopfc[\opt][\delay](\lambda)$	 without delay and with $\delay\in\{1,2,3\}$ ($r = 1$).}
	\label{fig:kstar-freq}
\end{figure}

\Cref{fig:kstar-freq} shows the optimal feedback gains without delays and with several values of the delay $\delay$, with $r = 1$.
The optimal gain $\Kopfc[\opt][\delay](\lambda)$ is smaller than the corresponding delay-free optimal gain and decreases with $\delay$,
hence the solution of~\eqref{eq:prob-optimal-control-decoupled-explicit} requires smaller control effort with long delays.
This agrees with the asymptotic results in \Crefns{prop:optimal-gain-tails,thm:optimal-gain-expensive-regime,prop:optimal-gain-small-delay},
suggesting that feedback becomes less relevant in the presence of delays.
The domain of the curves progressively shrinks leftwards because of the stability condition $\Aopf\delay<1$.

\section{Optimal Control of Spatially Invariant Systems with Communication Delays}\label{sec:centralized}

In this section, we provide our main results on optimal control of spatially invariant systems, such as those described in \cref{ex:reaction-diffusion,ex:circulant}, in the presence of communication delays that enter the feedback loop.
The analysis presented in \autoref{sec:analysis} -- characterizing optimal feedback gains for the \red{scalar} problem~\eqref{eq:prob-optimal-control-decoupled} -- serves as a foundation for the results on spatially invariant systems. 
The connection comes through the fact that the spatial Fourier transform decouples spatially invariant systems to a family of independent subsystems in the spatial frequency domain.
A spatially invariant system can be indeed thought of as a family of systems parametrized by the spatial frequency $\lambda$.
We consider the spatially invariant system
\begin{equation}\label{eq:model-spatially-invariant}
	\begin{aligned}
		\ddt{\psi}{x,t} &= \A{\psi} + \u{x}{t} + \noise{x}{t}\\
		\u{x}{t}	&=  -\Kc{\psi}{T}\\
		\z{x}{t} &=\begin{bmatrix}
			\state{x}{t}\\
			\sqrt{r} \, \u{x}{t}
		\end{bmatrix}
	\end{aligned}
\end{equation}
where $\psi(\cdot,t)\in L^2(\spacedom)$ is a scalar-valued function on a suitable space domain $\spacedom$ (see \cref{app:preliminaries}) and $\noise{\cdot}{t} \in L^2(\spacedom)$.
To ensure stabilizability of system~\eqref{eq:model-spatially-invariant},
we assume that the open-loop spatially invariant operator $\Aop:\mathcal{D}(\Aop) \subseteq L^2(\spacedom)\to L^2(\spacedom)$ is the infinitesimal generator of a $\mathcal{C}_0$-semigroup with continuous Fourier symbol $\Aopftot$ satisfying \cref{ass:A-even}~\cite{9683472},
and that the feedback controller $\Kop{}$ is a spatial convolution operator \revision{whose kernel $\Kopk:\spacedom\to\Real{}$ is an even function} as per \cref{ass:A-even}.
This,
because the spectra of stable systems of class~\eqref{eq:model-spatially-invariant-fourier} are strictly bounded away from zero~\cite{datko1978procedure},
ensures that the closed-loop system~\eqref{eq:model-spatially-invariant} is (exponentially) stable and makes the $\mathcal{H}_2$-norm from $v$ to $z$ well defined.
Denoting the latter by $J$,
we address the following optimal control design problem:
\begin{argmini!}
	{\Kop{}}
	{J(\Kop{})\protect\label{eq:prob-optimal-control-obj}}
	{\label{eq:prob-optimal-control}}
	{\Kop{\delay}^\opt \in }
	\addConstraint{\eqref{eq:model-spatially-invariant} \mbox{ is stable.}\protect\label{eq:prob-optimal-control-constr-dyn}}
\end{argmini!}
In the rest of this section,
we characterize the optimal feedback operator $\mathcal{K}_T^{\text{opt}}$,
which by construction is a spatial convolution,
by analyzing its kernel $K_T^{\text{opt}}(\cdot)$.
Specifically, we consider the expensive control regime ($r \to \infty$)  in \autoref{sec:expensive-control-spatial-structure} and the small communication delay regime ($T \to 0^+$) in \autoref{sec:small-delays-spatially-invariant}.

\subsection{Expensive Control Regime}\label{sec:expensive-control-spatial-structure}

For stable open-loop systems,
we characterize the optimal control kernel $\Kopk[\opt][\delay](\cdot)$ at the limit for $r\rightarrow\infty$,
\ie when state deviations are negligibly penalized compared to control effort.

\begin{thm}[Optimal controller in the expensive control regime]\label{thm:expensive-regime}
	Let $\Aop$ with Fourier symbol strictly bounded away from zero,
	\ie there exists $A_\text{th}<0$ such that $\Aopftot\le A_\text{th} \; \forall\lambda$.
	Define
	\begin{equation}\label{eq:kstar-expensive}
		\Kopk[\text{ex}][\delay] \doteq \sqrt{2\pi}\lr\Kopk[\text{ex}][0] \star g_\delay \rr
	\end{equation}
	where
	\begin{equation}\label{eq:kstar-expensive-no-del}
		\Kopk[\text{ex}][0] \doteq \ifourier{\dfrac{1}{2r|\Aopftot|}} \qquad g_\delay \doteq \ifourier{\e^{\delay\Aopf}}.
	\end{equation}
	Then,
	it holds
	\begin{equation}
		\lim_{r\rightarrow+\infty} \norm[L_2]{\Kopk[\opt][\delay] - \Kopk[\text{ex}][\delay]} = 0.
	\end{equation}
	\begin{proof}[Sketch of Proof]
		It follows from decoupling of~\eqref{eq:prob-optimal-control} to~\eqref{eq:prob-optimal-control-decoupled},
		\cref{thm:optimal-gain-expensive-regime},
		and properties of the Fourier transform.
		The hypothesis that $\Aopftot$ is strictly bounded away from zero ensures that \cref{thm:optimal-gain-expensive-regime} holds true at all spatial frequencies $\lambda$.
	\end{proof}
\end{thm}

The kernel $\Kopk[\text{ex}][0]$ in~\eqref{eq:kstar-expensive-no-del} approximates the optimal controller kernel in the expensive regime for the delay-free case,
according to~\eqref{eq:kstar-no-delay}.
\cref{thm:expensive-regime} establishes a relationship between the optimal controller with communication delays and the optimal controller in the absence of delays.
In fact,
expression~\eqref{eq:kstar-expensive} shows that the optimal controller kernel in the expensive control regime is computed as the cascade of a delay-aware filter followed by the optimal kernel for the delay-free case, \ie
\begin{equation}\label{eq:kstar-expensive-approx-convolution}
	\begin{aligned}
		\u{x}{t}	&= -\dfrac{1}{\sqrt{2\pi}}\lr\Kopk[\text{ex}][\delay] \star \psi\rr(x,t-\delay) \\
		&= -\lr\Kopk[\text{ex}][0] \star g_\delay \star \psi\rr(x,t-\delay).
	\end{aligned}
\end{equation}
This decomposition provides insight into how communication delays affect the spatial structure of the optimal feedback kernel.
The delay-aware controller uses the filter $g_\delay$ to counterbalance delays in the feedback loop.
If the spatial spread of $g_\delay$ is wider than that of the delay-free control kernel $\Kopk[\text{ex}][0]$,
truncation of the controller gains based on the spatial locality of $\Kopk[\text{ex}][0]$ may degrade performance, as it discards measurements relevant to handle delays.
In this case,
a distributed controller implementation based on truncation should mainly consider the spatial ``width'' of the delay-aware filter $g_\delay$.
A concrete visualization for the system considered in \cref{ex:reaction-diffusion} is discussed in~\autoref{sec:diffusion}.

Also,
the expression of $\Kopfc[\text{ex}][\delay]$ in~\eqref{eq:kstar-expensive-approx} shows that the controller kernel $\Kopk[\text{ex}][\delay]$ exhibits the same asymptotic spatial exponential decay of the delay-free controller (under mild assumptions on the open-loop dynamics).
This readily follows from~\cite{1017552} because $\Kopfc[\text{ex}][\delay]$ is the pointwise product of $\Kopfc[\text{ex}][0]$ with an exponential function of $\Aopf$,
which is analytic on the complex plane.

\subsubsection*{Performance gap}\label{sec:expensive-control-cost} 

The spatial frequency decoupling lets us compute the optimal $\mathcal{H}_2$-norm of problem~\eqref{eq:prob-optimal-control}.
The domain of $\lambda$ is $\Lambda$ and characterized in \cref{app:preliminaries}.
By continuity,
it holds:
\begin{equation}\label{eq:jstar-expensive}
	\begin{aligned}
		J\lr\Kop{\delay}^\text{ex}\rr	&= \int_{\Lambda}J_\lambda\lr\Kopfc[\text{ex}][\delay](\lambda)\rr\de\lambda\\
		&\sim\int_{\Lambda}\dfrac{1}{2|\Aopftot|} \lr1-\dfrac{\e^{2\delay\Aopftot}}{4r\Aopftot[2]}\rr\de\lambda.
	\end{aligned}
\end{equation}
The detailed derivation of~\eqref{eq:jstar-expensive} is deferred to~\cref{app:optimal-cost-expensive-regime}.
We compare the cost~\eqref{eq:jstar-expensive} with the $\mathcal{H}_2$-norm in the delay-free scenario under expensive regime:
\begin{equation}\label{eq:jstar-expensive-no-delay}
	J\lr\Kop{0}^\text{ex}\rr \sim \int_{\Lambda}\dfrac{1}{2|\Aopftot|} \lr1-\dfrac{1}{4r\Aopftot[2]}\rr\de\lambda.
\end{equation}
The correction factor $\e^{\delay\Aopftot}$ that accounts for delays smoothly inflates the $\mathcal{H}_2$-norm from~\eqref{eq:jstar-expensive-no-delay} to~\eqref{eq:jstar-expensive}.
The optimal cost under communication delays is larger by the following amount:
\begin{equation}\label{eq:jstar-diff}
	J\lr\Kop{\delay}^\text{ex}\rr - J\lr\Kop{0}^\text{ex}\rr
	\sim \Delta J^\text{ex}(\delay) \doteq \int_{\Lambda} \dfrac{1-\e^{2\delay\Aopftot}}{8r|\Aopftot[3]|}\de\lambda.
\end{equation}
The cost gap $\Delta J^\text{ex}(\delay)$ has the same order of magnitude of $r^{-1}$ irrespectively of the delay $\delay$.
This is explained because the delay is introduced through the control.
Because the control kernel $\Kopk[\text{ex}][\delay]$ is proportional to $r^{-1}$,
the cost gap has the same order of magnitude.
At the limit,
it holds
\begin{subequations}
	\begin{align}
		\lim_{\delay\rightarrow0^+} \Delta J^\text{ex}(\delay) &= 0\label{eq:cost-gap-T-0}\\
		\lim_{\delay\rightarrow+\infty} \Delta J^\text{ex}(\delay) &= \int_{\Lambda}\dfrac{1}{8r|\Aopftot[3]|}\de\lambda.\label{eq:cost-gap-T-infty}
	\end{align}
\end{subequations}
Limit~\eqref{eq:cost-gap-T-0} shows that,
as the delay vanishes,
the gap between the optimal $\mathcal{H}_2$-norms vanishes as well because the optimal controller recovers the delay-free one.
On the other hand,
as $\delay$ grows,
the optimal control becomes negligible and the limit~\eqref{eq:cost-gap-T-infty} yields the $\mathcal{H}_2$-norm of the stable open-loop system.

\subsection{Small Communication Delays}\label{sec:small-delays-spatially-invariant}

We now consider the case where communication delays affecting the feedback loop are negligible compared to the timescale of the open-loop dynamics.

\begin{thm}[Optimal control kernel in the small delay regime]\label{thm:small-delays}
	Define
	\begin{equation}\label{eq:kstar-small-T-space}
		\Kopk[\text{sd}][\delay] \doteq \ls(\mathcal{I} - \Aop\delay)\Kopk[\opt][0]\rs - \dfrac{\delay}{r}\delta
	\end{equation}
	where $\mathcal{I}$ is the identity operator,
	$\Kopk[\opt][0]$ is the optimal control kernel for problem~\eqref{eq:prob-optimal-control} with $\delay=0$,
	and $\delta$ is the Dirac delta.
	If $|\Aopftot|<\Aopf[\text{th}] \; \forall\lambda$ for some $\Aopf[\text{th}]\in\Real{}$,
	it holds
	\begin{equation}
		\lim_{\delay\rightarrow0^+} \norm[L_2]{\Kopk[\opt][\delay] - \Kopk[\text{sd}][\delay]} = 0.
	\end{equation}
	\begin{proof}[Sketch of Proof]
		It follows from decoupling of~\eqref{eq:prob-optimal-control} to~\eqref{eq:prob-optimal-control-decoupled},
		\cref{prop:optimal-gain-small-delay},
		and properties of the Fourier transform.
		The hypothesis that $\Aopftot$ is uniformly bounded ensures that \cref{prop:optimal-gain-small-delay} holds true at all spatial frequencies $\lambda$.
	\end{proof}
\end{thm}

Expression~\eqref{eq:kstar-small-T-space} shows that,
under small delays,
the optimal delay-free control kernel is perturbed by a $\delay$-scaled version of the open-loop operator $\Aop$ and by an additive negative Dirac delta.
The bounded Fourier symbol $\Aopf$ ensures that the delay $\delay$ is small compared to the characteristic timescale of the open-loop dynamics,
corresponding to the regime $\alpha_*\ll1$ discussed in \autoref{sec:tails}.
Further,
the scaled operator $\Aop\delay$ is bounded,
so that the operator $\mathcal{I} - \Aop\delay$ in~\eqref{eq:kstar-small-T-space} actually produces a ``small'' perturbation of the delay-free kernel $\Kopk[\text{op}][0]$ which varies linearly with $\delay$.

Also,
expression~\eqref{eq:kstar-small-T-approx} of $\Kopfc[\text{sd}][\delay]$ in spatial frequency domain is affine with respect to the spatial Fourier transform of the delay-free optimal controller.
This allows us to determine the asymptotic exponential decay of the control kernel $\Kopk[\text{sd}][\delay]$,
which readily follows from the result in~\cite{1017552} similarly to what observed for the expensive control regime.

\begin{rem}[Continuous- vs. discrete-time]
	\revision{\Cref{thm:expensive-regime,thm:small-delays} build on analysis of scalar systems in \autoref{sec:analysis} which uses the closed-form expression of the $\mathcal{H}_2$-norm $J_\lambda$ derived in~\cite{KuchlerLangevinEqs}.
		This holds for continuous-time dynamics and prevents us from readily extending our results to discrete time.}
\end{rem}

\section{Example One: Reaction-Diffusion Process}\label{sec:diffusion}

To make the analysis in~\autoref{sec:centralized} more concrete,
we now consider the spatially invariant system in \cref{ex:reaction-diffusion},
namely a real scalar reaction-diffusion process over the real line.

\subsection{Optimal Control Kernel in Expensive Control Regime}\label{sec:diffusion-expensive}

\begin{figure}
	\centering
	\begin{subfigure}{.5\linewidth}
		\centering
		\includegraphics[width=\linewidth]{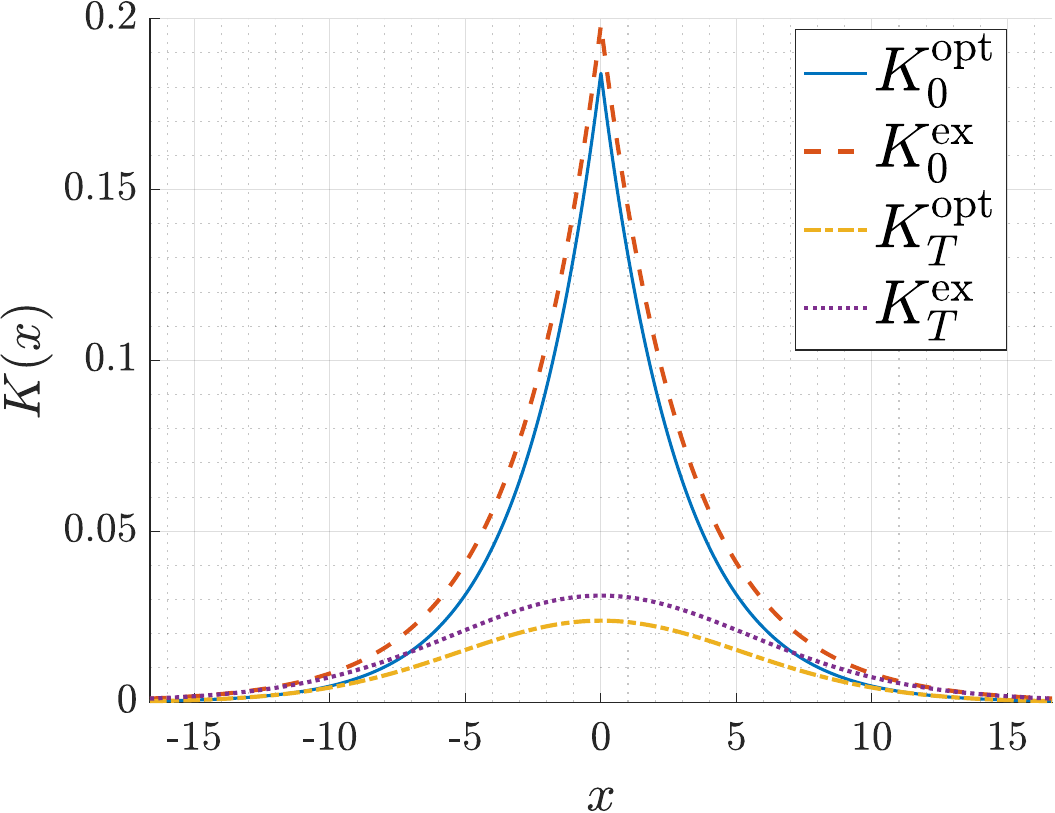}
		\caption{Control weight $r=1$.}
		\label{fig:kstar_space_d_10_c_1_T_1_r_1}
	\end{subfigure}%
	\hfill
	\begin{subfigure}{.5\linewidth}
		\centering		
		\includegraphics[width=\linewidth]{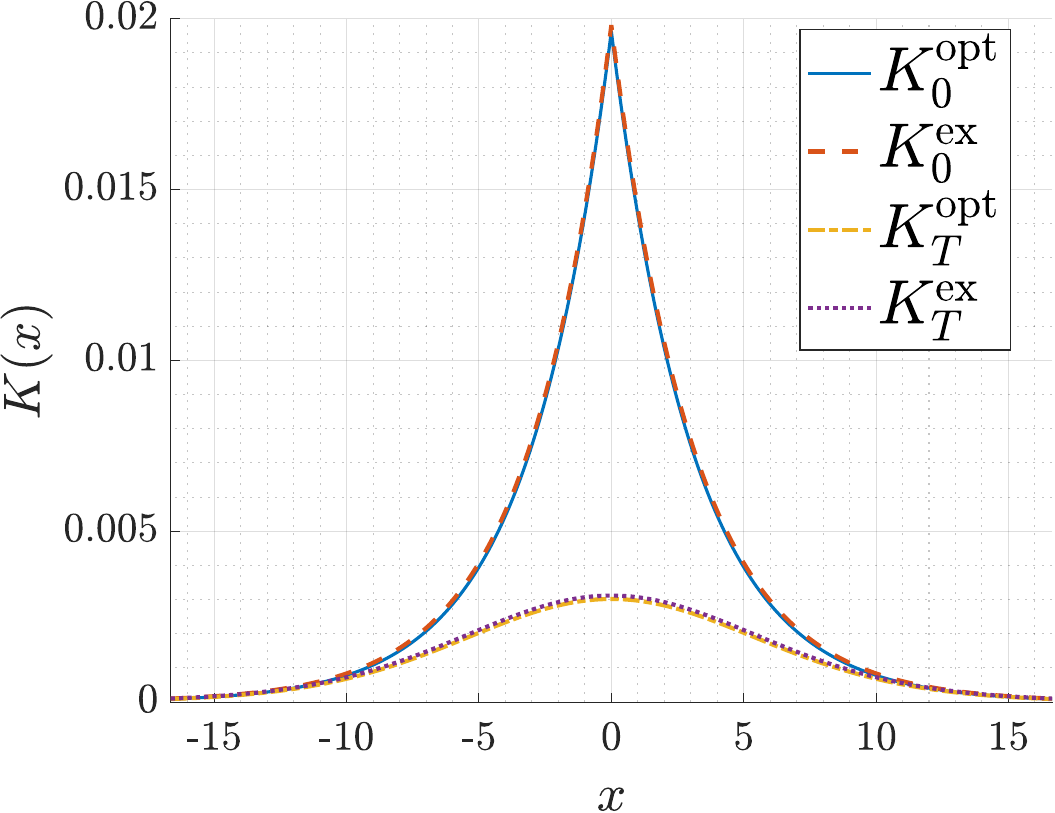}
		\caption{Control weight $r=10$.}
		\label{fig:kstar_space_d_10_c_1_T_1_r_10}
	\end{subfigure}
	\caption{Optimal controllers without delay vs. delay $\delay=1$ and their approximations in the expensive-control regime with $d=10$ and $c=1$.}
	\label{fig:kstar-diffusion}
\end{figure}

In the expensive control regime,
the optimal controller is given by~\cref{thm:expensive-regime} as $\Kopk[\text{ex}][\delay] = \sqrt{2\pi} \lr\Kopk[\text{ex}][0]\star g_\delay\rr$,
where the optimal controller kernel for the delay-free case is
\begin{equation}\label{eq:kstar-expensive-diffusion-no-delay-space}
	\Kopk[\text{ex}][0](x) = \ifourier[x]{\dfrac{1}{2r\lr d\lambda^2+c\rr}} = \dfrac{1}{2r}\sqrt{\dfrac{\pi}{2dc}}\e^{-\sqrt{\frac{c}{d}}|x|}
\end{equation}
and the delay-aware filter is
\begin{equation}\label{eq:kstar-expensive-diffusion-gaussian-kernel}
	g_\delay(x) = \ifourier[x]{\e^{-\delay\lr d\lambda^2+c\rr}} = \dfrac{\e^{-c\delay}}{\sqrt{2d\delay}}\e^{-\frac{x^2}{4d\delay}}.
\end{equation}
Tedious but straightforward computations yield
\begin{equation}\label{eq:kstar-expensive-diffusion-space}
	\Kopk[\text{ex}][\delay](x) = \dfrac{1}{2r}\sqrt{\dfrac{\pi}{2dc}}\lr\phi(x) + \phi(-x)\rr
\end{equation}
where we define
\begin{equation}\label{eq:kstar-expensive-diffusion-space-gamma}
	\phi(x) \doteq \dfrac{\e^{\sqrt{\frac{c}{d}}x}}{2}\lr1+\erf\lr-\dfrac{x}{2\sqrt{d\delay}}-\sqrt{c\delay}\rr\rr.
\end{equation}
\revision[According to \cref{thm:expensive-regime},
function {$\Kopk[\text{ex}][\delay]$ in~\eqref{eq:kstar-expensive-diffusion-space} asymptotically approaches the optimal controller kernel $\Kopk[\opt][\delay]$ as $r$ increases.}]{}

\subsection{Numerical Solutions and Approximations for Design}\label{sec:diffusion-numerical-and-design}

We now numerically compute the optimal convolution kernel $\Kopk[\opt][\delay]$ from the solution to the scalar problem~\eqref{eq:prob-optimal-control-decoupled} for system~\eqref{eq:model-diffusion-ex}--\eqref{eq:controller-convolution-ex},
and compare it with the analytical optimal kernel for expensive control regime $\Kopk[\text{ex}][\delay]$ in~\eqref{eq:kstar-expensive-diffusion-space}.
The kernel $\Kopk[\text{ex}][\delay]$ is plotted in~\autoref{fig:kstar-diffusion} (purple dotted)
together with the optimal kernel $\Kopk[\opt][\delay]$ (yellow dashed-dotted) for delay $\delay=1$.
Also,
we compare the delay-free optimal control kernel $\Kopk[\opt][0]$ obtained from~\eqref{eq:kstar-no-delay} (blue solid)
with the analytical delay-free kernel for expensive regime $\Kopk[\text{ex}][0]$ in~\eqref{eq:kstar-expensive-diffusion-no-delay-space} (red dashed).
Both approximations for expensive control regime are close to the optimal kernels already for $r = 10$.
While the feedback gains enjoy the same asymptotic exponential decay,
communication delays induce fundamentally different shapes \revision[of the two controller kernels]{} about the origin.

We next study expression~\eqref{eq:kstar-expensive-diffusion-space} to achieve deeper analytical insight about spatial locality of the optimal control kernel.

\subsubsection{Approximation about origin}\label{sec:diffusion-about-origin}

Because $\Kopk[\text{ex}][\delay]\in\mathcal{C}^\infty(\Real{})$,
we approximate it as a parabola via a standard Taylor expansion about $x=0$.
We formalize this approximation as follows.

\begin{lemma}[Approximation about the origin]\label{lem:kstar-expensive-diffusion-parabola}
	Let
	\begin{equation}\label{eq:kstar-expensive-approximation-diffusion-derivative-coefs}
		\begin{aligned}
			D_0	&\doteq 1 - \erf\lr\sqrt{c\delay}\rr\\
			D_2 &\doteq \dfrac{cD_0}{2d} - \sqrt{\dfrac{c}{\pi\delay}}\dfrac{\e^{-c\delay}}{2d}.
		\end{aligned}
	\end{equation}
	Then,
	it holds
	\begin{equation}\label{eq:kstar-expensive-approximation-origin}
		\Kopk[\text{ex}][\delay](x) = \Kopk[\text{ex}][0](0) \lr D_0 +D_2x^2\rr + o(x^3) \quad \forall\,x\in\Real{}.
	\end{equation}
	\begin{proof}[Sketch of proof]
		Because $\Kopk[\text{ex}][\delay]$ is even,
		all derivatives with odd order are zero at $x=0$.
		Computing the second-order Taylor expansion about $x=0$ yields expression~\eqref{eq:kstar-expensive-approximation-origin}.
	\end{proof}
\end{lemma}

\cref{lem:kstar-expensive-diffusion-parabola} shows that the optimal kernel in the expensive control regime is approximated by a concave quadratic function about the origin.
The relative gap between the optimal feedback gains for delay-free and delay systems at $x=0$ (corresponding to the location where the control is applied)
amounts to
\begin{equation}
	\Delta\Kopk[\text{ex}][\delay](0) \doteq \dfrac{\Kopk[\text{ex}][0](0) - \Kopk[\text{ex}][\delay](0)}{\Kopk[\text{ex}][0](0)} = \erf\lr\sqrt{c\delay}\rr.
\end{equation}
The gap $\Delta\Kopk[\text{ex}][\delay](0)$ increases with $c$ and $\delay$.
This suggests that,
in the presence of  communication delays,
the optimal feedback gain at the actuator location become smaller as the delay and the reaction coefficient increase,
\revision[requiring less aggressive control]{} compared to the delay-free scenario.
Notably,
a more ``intense'' reaction term $-c\state{x}{t}$ in the dynamics~\eqref{eq:model-diffusion-ex} induces more conservative control actions
for every $\delay>0$.
This might be explained because this term stabilizes the system,
hence a larger reaction coefficient allows to reduce the control cost while also inducing a milder effect of delays on the closed-loop dynamics.
In contrast,
the diffusion term $ \frac{\partial^2\state{x}{t}}{\partial x^2}$ does not play a role in the gap $\Delta\Kopk[\text{ex}][\delay](0)$.
An explanation is that $\Kopk[\text{ex}][\delay](0)$ multiplies the state at the same location $x$ where the control $\u{x}{t}$ is applied,
\ie the product $\Kopk[\text{ex}][\delay](0)\state{x}{t-\delay}$ appears in~\eqref{eq:controller-convolution-ex},
whereas the diffusion term generates dynamic couplings among all spatial locations of the state $\state{\cdot}{t}$.

\cref{lem:kstar-expensive-diffusion-parabola} also reveals that the (negative) curvature $D_2$ increases with $d$ and $\delay$.
The control kernel becomes ``flatter'' about the origin with both longer communication delays and larger diffusion coefficient.
This suggests that the controller needs to incorporate state measurements from far away if the process is very sensitive to state variations across space (large $d$) or if the feedback measurements incur long delays $\delay$.

Overall,
the behavior of $D_0$ and $D_2$ that increase with $\delay$ means that the control kernel becomes flatter (shorter and more spread out) about the origin as the delay increases;
\revision[this requires more resources to communicate with state locations farther away and reduces the impact of feedback,
making control actions less aggressive and performing than delay-free control]{}.

\subsubsection{Approximation of tails}\label{sec:diffusion-tails}

As noted in \autoref{sec:expensive-control-spatial-structure},
the asymptotic spatial decay of the control kernel $\Kopk[\text{ex}][\delay]$ is (at least) exponential like the delay-free case.
In fact,
evaluating the limit of expression~\eqref{eq:kstar-expensive-diffusion-space} reveals that the kernel $\Kopk[\text{ex}][\delay]$ asymptotically approaches the delay-free kernel $\Kopk[\text{ex}][0]$ as $|x|\rightarrow+\infty$.
The following result quantifies how the delay-free optimal kernel approximates the delay-aware one for finite values of $x$.

\begin{lemma}[Approximation of tails]\label{lem:kstar-expensive-diffusion-tails}
	It holds
	\begin{equation}\label{eq:kstar-expensive-diffusion-approximation-x-large}
		\Kopk[\text{ex}][\delay](x) = \Kopk[\text{ex}][0](x)(1 + R(x)) \quad \forall\,x\in\Real{}
	\end{equation}
	where the function $R$ is negative and lower bounded as
	\begin{multline}\label{eq:kstar-expensive-diffusion-approximation-x-large-remainder-bound}
		R(x) > \dfrac{\e^{-\frac{1}{2}\lr\frac{x^2}{2d\delay}+2c\delay\rr}}{\sqrt{2\pi}}\lr\dfrac{\sqrt{2d\delay}}{x+2\sqrt{dc}\delay}-\dfrac{\sqrt{2d\delay}^3}{(x+2\sqrt{dc}\delay)^3}\rr\\
		-\dfrac{\e^{-\frac{1}{2}\lr\frac{x}{\sqrt{2d\delay}}-\sqrt{2c\delay}\rr^2}}{\sqrt{2\pi}}\dfrac{\sqrt{2d\delay}}{x-2\sqrt{dc}\delay}
	\end{multline}
	and has limit
	\begin{equation}\label{eq:kstar-expensive-diffusion-approximation-x-large-remainder-limit}
		\lim_{|x|\rightarrow+\infty}R(x) = 0.
	\end{equation}
	\begin{proof}[Sketch of proof]
		It follows by applying the divergent asymptotic approximation of $\Phi(x)=\nicefrac{1}{2}\lr1+\erf(x)\rr$ to $\nicefrac{\Kopk[\text{ex}][\delay](x)}{\Kopk[\text{ex}][0](x)}$ and performing straightforward algebraic manipulations.
	\end{proof}
\end{lemma}

\cref{lem:kstar-expensive-diffusion-tails} shows that the convolution kernel of the delay-aware controller can be effectively approximated by the standard delay-free kernel for large values of $x$.
The accuracy of this approximation is quantified by the function $R$ that is lower bounded in~\eqref{eq:kstar-expensive-diffusion-approximation-x-large-remainder-bound}.
The exponential functions in~\eqref{eq:kstar-expensive-diffusion-approximation-x-large-remainder-bound} derive from function $\phi$ in~\eqref{eq:kstar-expensive-diffusion-space-gamma}.
For $|x|$ large enough,
the exponents in~\eqref{eq:kstar-expensive-diffusion-approximation-x-large-remainder-bound} decrease with the delay through the product $c\delay$,
suggesting that the combination of reaction term and delay has the greatest impact in the asymptotic spatial decay of the kernel.

\subsubsection{Practical design guidelines}\label{sec:diffusion-design}

The asymptotic approximations given by~\cref{lem:kstar-expensive-diffusion-parabola,lem:kstar-expensive-diffusion-tails}
lend themselves to practical guidelines for delay-aware controller design by finding intervals where the approximations are sufficiently accurate.
A possible rule of thumb is given by imposing that the quadratic function of the Taylor-based approximation~\eqref{eq:kstar-expensive-approximation-origin}
dominates the higher-order terms,
whereas the residual term $R$ in the approximation~\eqref{eq:kstar-expensive-diffusion-approximation-x-large}
can be neglected by imposing that the exponents of the two exponential functions are small enough.
More precisely,
let us define
\begin{align}
	D_4 			&\doteq 2\dfrac{c^2D_0}{d^2} - \sqrt{\dfrac{c}{\pi\delay}}\dfrac{\e^{-c\delay}}{d^2}\lr c+\dfrac{1}{2\delay}\rr\\
	x_\text{th,1}	&\doteq \sqrt{\dfrac{12}{|D_4|}\lr D_2 + \sqrt{D_2^2 + \dfrac{D_0|D_4|}{6}}\rr}\label{eq:x-th1}\\
	x_\text{th,2}	&\doteq2\lr\sqrt{d\delay} + \sqrt{cd}\delay\rr.\label{eq:x-th2}
\end{align}
The thresholds $x_\text{th,1}$ and $x_\text{th,2}$ are calculated such that the asymptotic expressions in~\cref{lem:kstar-expensive-diffusion-parabola,lem:kstar-expensive-diffusion-tails} yield good approximations in the intervals $(0, x_\text{th,1})$ and $(x_\text{th,2},\infty)$,
respectively,
such that the remainders are negligible.
Then,
for $|x| <\alpha x_\text{th,1}$ and $|x| > \beta x_\text{th,2}$,
where $\alpha \in(0,1)$ and $\beta>0$ are design parameters such that $\alpha x_\text{th,1} \le \beta x_\text{th,2}$,
the control kernel $\Kopk[\text{ex}][\delay]$ can be approximated by the function $\widetilde{\Kopk[\text{ex}][\delay]}$ defined as follows:
\begin{equation}\label{eq:kstar-approx-design}
	\widetilde{\Kopk[\text{ex}][\delay]}(x) = 
	\begin{cases}
		\Kopk[\text{ex}][0](0) \lr D_0 +D_2x^2\rr	& |x| \le \alpha x_\text{th,1}\\
		\Kopk[\text{ex}][0](x)						& |x| \ge \beta x_\text{th,2}.
	\end{cases}
\end{equation}
The parameters $\alpha$ and $\beta$ are meant to trade a simple design (large $\alpha$ and small $\beta$) for an accurate approximation (small $\alpha$ and large $\beta$).
The feedback gain in the interval $(\alpha x_\text{th,1}, \beta x_\text{th,2})$ may be chosen in practice by interpolating the two cases in~\eqref{eq:kstar-approx-design}.

\begin{figure}
	\centering
	\begin{subfigure}{0.5\linewidth}
		\centering
		\includegraphics[width=\linewidth]{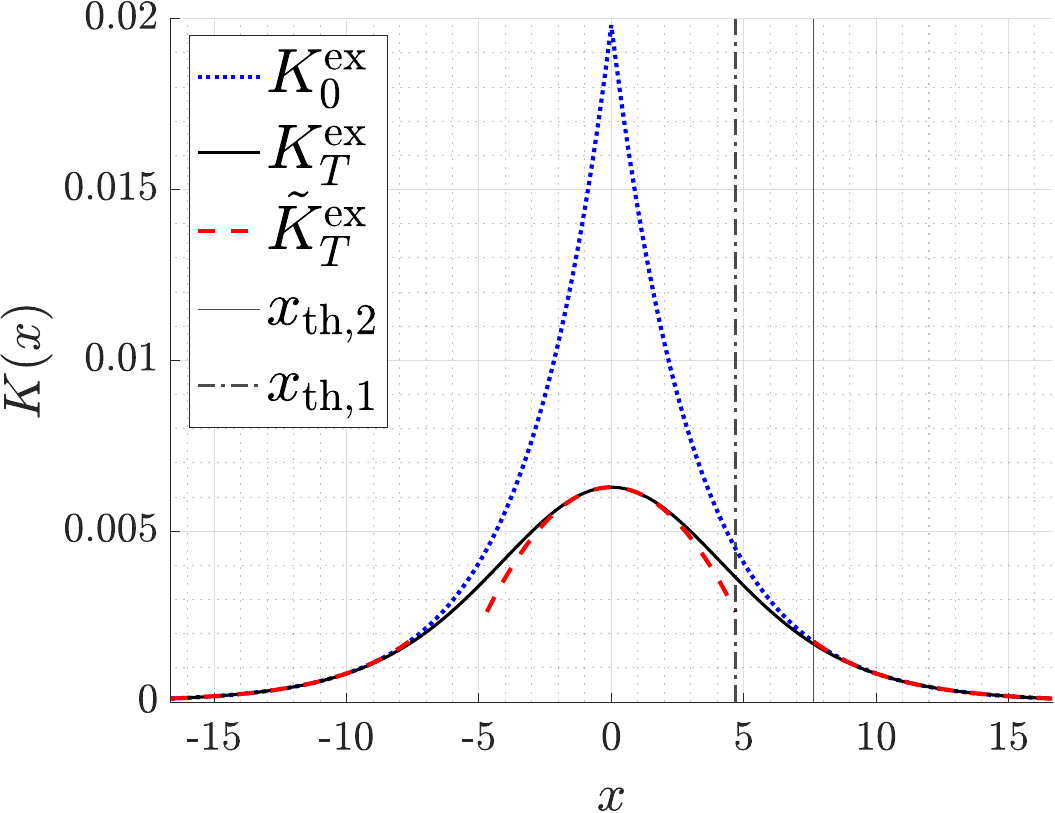}
		\caption{Delay $\delay=0.5$.}
		\label{fig:kstar-diffusion-expensive-approximation-for-design-0.5}
	\end{subfigure}%
	\hfill
	\begin{subfigure}{0.5\linewidth}
		\centering
		\includegraphics[width=\linewidth]{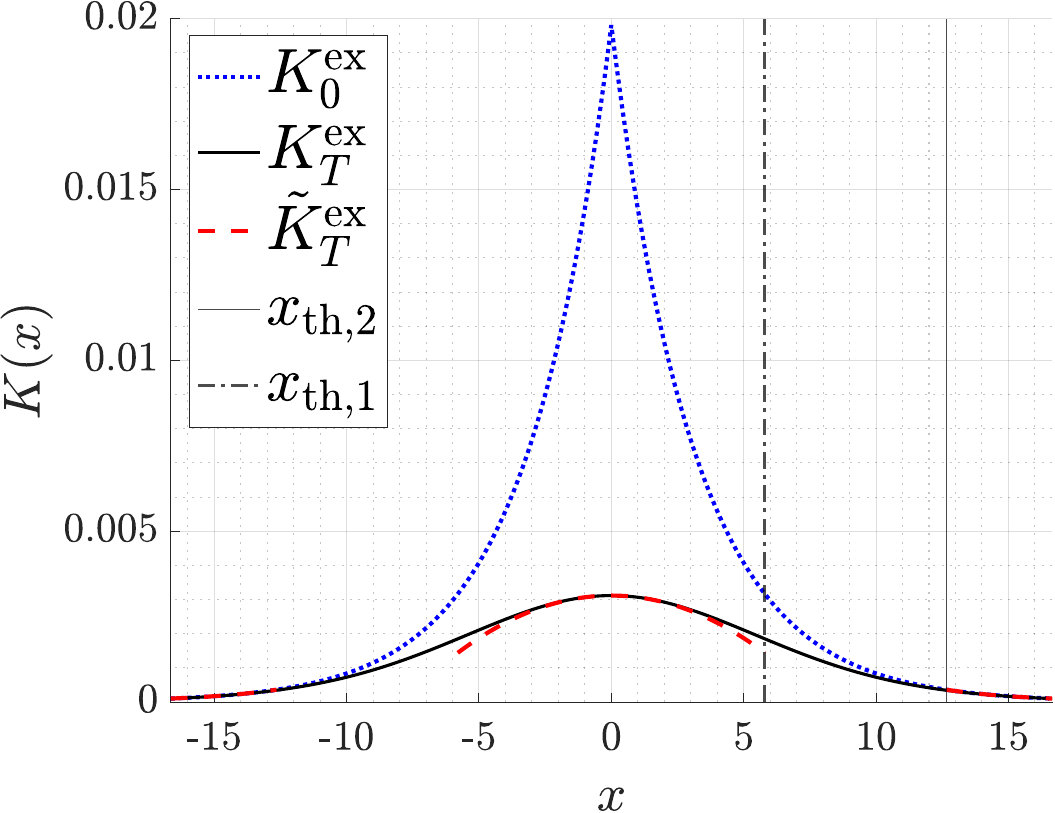}
		\caption{Delay $\delay=1$.}
		\label{fig:kstar-diffusion-expensive-approximation-for-design-1}
	\end{subfigure}
	\caption{Controllers in the expensive-control regime vs. design approximation $\widetilde{\Kopk[\text{ex}][\delay]}(x)$ with $d=10$, $c=1$, and $r=10$.}
	\label{fig:kstar-diffusion-expensive-approximation-for-design}
\end{figure}

\Cref{fig:kstar-diffusion-expensive-approximation-for-design} compares the delay-free controller $\Kopk[\text{ex}][0]$ (dotted blue),
the delay-aware controller $\Kopk[\text{ex}][\delay]$ (solid black),
and the approximation $\widetilde{\Kopk[\text{ex}][\delay]}$ (dashed red) for $\alpha=\beta=1$, together with the thresholds $x_{\text{th},1}$ and $x_{\text{th},2}$.
The two cases of $\widetilde{\Kopk[\text{ex}][\delay]}$ in~\eqref{eq:kstar-approx-design} closely approach the optimal kernel $\Kopk[\text{ex}][\delay]$ as $x$ approaches the origin (first case) or grows large in magnitude (second case).

\subsubsection{Truncation of optimal control kernel}\label{sec:diffusion-truncation}

The exponentially decaying control kernel $\Kopk[\text{ex}][\delay]$ suggests that,
as pointed out in~\cite{1017552},
a simple yet effective way to implement a distributed controller architecture is spatial truncation of the feedback kernel gains.
However,
this has to be done with care because such a truncation can degrade control performance,
and can cause instability of the closed-loop system in general~\cite{4623272}.

In~\cite{1017552,9683472},
it is shown that performance loss for the delay-free scenario diminishes asymptotically exponentially with the truncation width as a consequence of the control kernel decay and small-gain arguments.
In expensive regime,
the kernel $\Kopk[\text{ex}][0]$ in~\eqref{eq:kstar-expensive-diffusion-no-delay-space} decays exponentially in space with rate $\sqrt{\nicefrac{c}{d}}$.
Hence,
a practical truncation rule that does not account for delays may discard spatial locations farther than $x_{\text{th},0}\doteq \gamma\sqrt{\nicefrac{d}{c}}$,
where $\nicefrac{1}{\gamma}\le1$ reflects the performance loss.

However,
as shown in~\autoref{fig:kstar-diffusion},
communication delays play a crucial role in shaping the optimal controller kernel,
which calls for further attention in implementing spatial truncation.
While a precise evaluation needs to be done numerically,
an approximate yet insightful assessment is conveyed by properties of convolution.
As discussed in~\autoref{sec:expensive-control-spatial-structure},
the delay-aware controller $\Kop{\delay}^\text{ex}$ is equivalent to the cascade of the optimal controller without delays~\eqref{eq:kstar-expensive-diffusion-no-delay-space} and the delay-aware filter~\eqref{eq:kstar-expensive-diffusion-gaussian-kernel},
\begin{equation}
	\begin{aligned}
		\u{x}{t}	&= -[\Kopk[\text{ex}][0] \star g_\delay \star \psi](x,t-\delay)\\
		&\propto -\e^{-\sqrt{\frac{c}{d}}|x|} \star \e^{-\frac{x^2}{4d\delay}} \star \state{x}{t-\delay}.
	\end{aligned}
\end{equation}
The expression above highlights that spatial truncation of the kernel $\Kopk[\text{ex}][\delay]$ affects both filtering actions.
Loss in performance is caused by discarding state measurements that are relevant to the delay-free control through $\Kopk[\text{ex}][0]$ and/or to the delay-aware filter through $g_{\delay}$.
The delay $\delay$ plays a fundamental role in determining which of these two kernels is more demanding in terms of spatial width,
which quantifies the effect of truncation.
Similarly to the discussion above to truncate $\Kopk[\text{ex}][0]$,
a rule of thumb to truncate the Gaussian kernel $g_\delay$ considers that most of its area is within $x_{\text{th},\delay}\doteq \kappa\sqrt{2d\delay}$ from the origin,
where $\kappa$ is usually $2$ or $3$.
Then,
a possible way to determine how delays affect performance under truncation of the control kernel as compared to the delay-free case uses the inequality
\begin{equation}\label{eq:kstar-diffusion-expensive-truncation-ineq}
	\sqrt{2c\delay} > \dfrac{\gamma}{\kappa}.
\end{equation}
If~\eqref{eq:kstar-diffusion-expensive-truncation-ineq} holds for user-defined parameters $\gamma$ and $\kappa$,
the truncated delay-aware filter kernel $g_{\delay}$ is wider than the truncated delay-free controller kernel $\Kopk[\text{ex}][0]$,
\revision[meaning that handling feedback delays imposes to retain state measurements from farther away compared to delay-free control.
Thus,
in this case,
the kernel {$\Kopk[\text{ex}][\delay]$} should be truncated according to the truncation threshold $x_{\text{th},\delay}$ of the delay-aware filter $g_\delay$,
and not of the delay-free kernel {$\Kopk[\text{ex}][0]$},
to maintain the desired near-optimal performance.]{}

\begin{figure}
	\centering
	\begin{subfigure}{.5\linewidth}
		\centering
		\includegraphics[width=.95\linewidth]{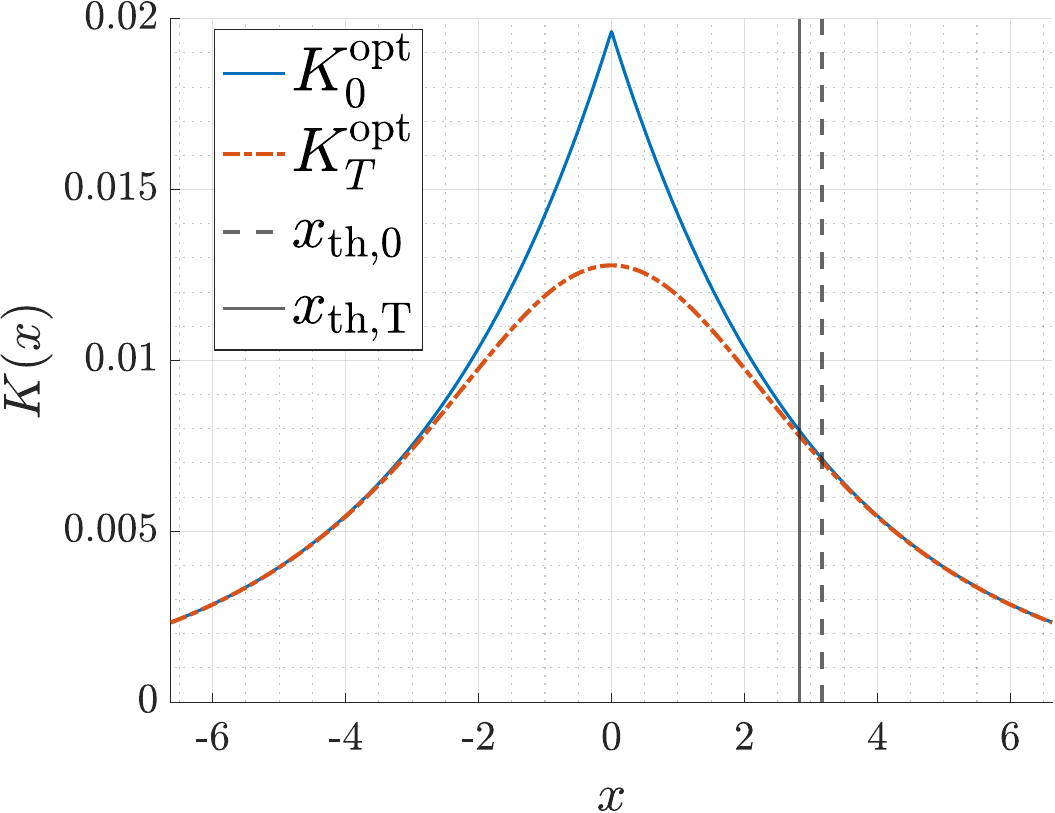}
		\caption{Delay $\delay=0.1$.}
		\label{fig:kstar_space_trunc_d_10_c_1_T_0.1_r_10_k_2_g_1}
	\end{subfigure}%
	\hfil
	\begin{subfigure}{.5\linewidth}
		\centering
		\includegraphics[width=.95\linewidth]
		{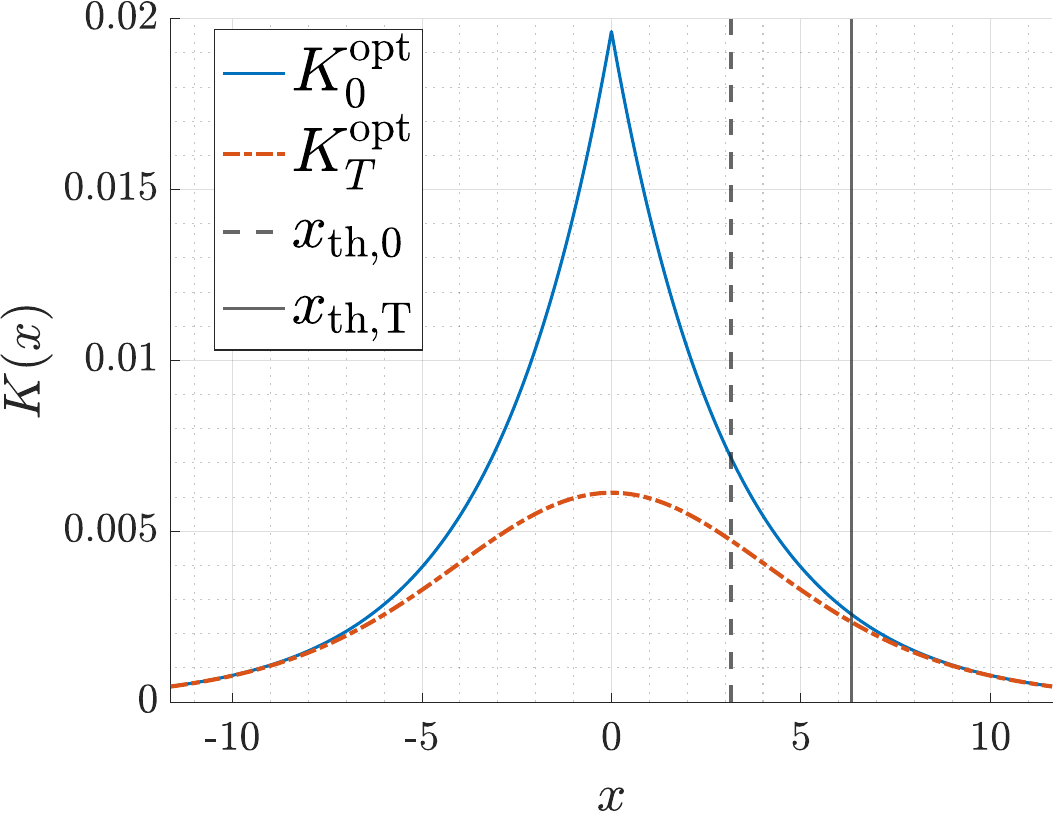}
		\caption{Delay $\delay=0.5$.}
		\label{fig:kstar_space_trunc_d_10_c_1_T_0.5_r_10_k_2_g_1}
	\end{subfigure}
	\caption{Optimal control kernels without delay vs. with nonzero delay together with truncation thresholds with $d=10$, $c=1$, $r=10$, $\kappa=2$, and $\gamma=1$.}
	\label{fig:kstar-truncation}
\end{figure}

Two opposite situations according to condition~\eqref{eq:kstar-diffusion-expensive-truncation-ineq} are illustrated in~\autoref{fig:kstar-truncation},
where the optimal delay-free kernel $\Kopk[\opt][0]$ and the optimal kernel under delays $\Kopk[\opt][\delay]$ are represented together with the two truncation thresholds $x_{\text{th},0}=\sqrt{\nicefrac{d}{c}}$,
which neglects communication delays,
and $x_{\text{th},\delay}=\sqrt{2d\delay}$,
which considers delays.
With $\delay=0.1$,
truncating the feedback gains without considering delays negligibly affects control performance because the delay-aware kernel $g_\delay$ is concentrated within the delay-free threshold $x_{\text{th},0}$; see \autoref{fig:kstar_space_trunc_d_10_c_1_T_0.1_r_10_k_2_g_1}.
However,
\autoref{fig:kstar_space_trunc_d_10_c_1_T_0.5_r_10_k_2_g_1} shows that the delay-aware truncation threshold $x_{\text{th},\delay}$ is about twice the delay-free threshold $x_{\text{th},0}$ if $\delay=0.5$,
thus the delay-free truncation would degrade performance because it discard state measurements relevant for delay-aware feedback.

The diffusion coefficient $d$ does not appear in the inequality~\eqref{eq:kstar-diffusion-expensive-truncation-ineq}.
Similarly to the discussion in \autoref{sec:diffusion-tails},
this suggests that the biggest contribution to the shape difference between delay-free and delay-aware control kernels is given by the spatially uncoupled reaction term,
while the diffusion term contributes similarly with and without delays.

\section{Example Two: Circulant Multi-Agent System}\label{sec:circulant}

We now turn to \cref{ex:circulant} that involves agents coupled through circulant matrices.
Given vectors $A\in\Real{N}$ and $\Kopk{}{}\in\Real{N}$,
we denote the open-loop and feedback gain matrices in~\eqref{eq:circulant-dynamics} by $\Aop = \mathrm{circ}(A)$ and $\Kop{} = \mathrm{circ}(\Kopk{}{})$,
respectively,
where $\mathrm{circ}(C)$ is the circulant matrix with the vector $C$ in the first row.
Thus,
the vectors $A$ and $K$ respectively represent open-loop couplings and control feedback gains between agents.
The Fourier symbol $\hat{C}_\lambda$ of the matrix $\mathrm{circ}(C)$ is the Discrete Fourier Transform of $C$~\cite{ballotta2023tcns}.
In view of \cref{ass:A-even},
we impose that $\Aop$ and $\Kop{}$ are symmetric,
and it follows that $A(i) = A(N-i)$ and $K(i) = K(N-i)$ for $i\in\{1,\dots,N-1\}$.
For agent $i$,
it holds
\begin{multline}
	\dfrac{\de\state{i}{t}}{\de t} = \sum_{j=0}^{N-1} A(|i-j|_N)\state{j}{t} \\
	-\sum_{j=0}^{N-1}  K(|i-j|_N)\state{j}{t-\delay} + \noise{i}{t}
\end{multline}
where $|c|_N \doteq |c| \mod \ceil{\nicefrac{(N-1)}{2}}$ represents distance $|c|$ in the ring topology,
and we respectively denote the first coordinates of $A$ and $\Kopk{}{}$ by $A(0)$ and $\Kopk{}{}(0)$,
corresponding to the autonomous dynamics of each agent.

\subsection*{Optimal Feedback Gain Matrix with Small Delays}\label{sec:circulant-small-delays}

We use \cref{thm:small-delays} to approximate the optimal feedback gain matrix $\Kop{\delay}^\opt$ when the communication delay is small.
Note that the hypothesis $|\Aopf|<\Aopf[\text{th}]$ is trivially satisfied by choosing $\Aopf[\text{th}] = \max_\lambda|\Aopf| + \epsilon$ for any $\epsilon>0$.
Directly applying~\eqref{eq:kstar-small-T-approx} yields the following approximation of the eigenvalues of $\mathcal{K}$:
\begin{equation}\label{eq:circulant-small-delay-eigs}
	\Kopfc[\opt][\delay](\lambda) \sim \Kopfc[\text{sd}][\delay](\lambda) = \lr 1 - \Aopf\delay\rr\Kopfc[\opt][0](\lambda) - \dfrac{\delay}{r}
\end{equation}
where $\Kopfc[\opt][0](\lambda), \lambda=1,\dots,N,$ are the eigenvalues of the optimal feedback gain matrix without delays $\Kop{0}^\opt$.
As observed in~\autoref{sec:small-delays},
the perturbed eigenvalues under delays $\Kopfc[\text{sd}][\delay]$ are smaller than $\Kopfc[\opt][0]$,
inducing weaker control actions at all spatial frequencies.
In particular,
large perturbations are caused by large positive eigenvalues $\Aopf$,
while negative eigenvalues $\Aopf$,
which correspond to stable open-loop subsystems in frequency domain,
partially mitigate the perturbation because $-\Aopf$ has sign opposite to $-\nicefrac{\delay}{r}$.

From~\eqref{eq:kstar-small-T-space},
we approximate the optimal feedback gains as
\begin{equation}\label{eq:circulant-small-delay-gains}
	\Kopk[\opt][\delay] \sim \Kopk[\text{sd}][\delay] = (\mathcal{I}_N - \mathcal{A}\delay)\Kopk[\opt][0] - \dfrac{\delay}{r}e_0
\end{equation}
where $\mathcal{I}_N$ is the identity matrix of dimension $N$
and  $e_0\in\Real{N}$ is the canonical vector with first coordinate equal to $1$ and zero elsewhere.
Equation~\eqref{eq:circulant-small-delay-gains} shows that the feedback gains associated with measurements received from neighboring agents are reduced through multiplication by matrix $\Aop$,
whereas the self-gain of each agent is further decreased by the term $\nicefrac{\delay}{r}$,
which increases with the delay and decreases with the control weight.
This suggests that longer delays cause smaller feedback gains,
while heavier penalization of the control effort induces smaller differences between the self-gains of the delay-aware and delay-free optimal feedback gain matrix.

\begin{figure}
	\centering
	\begin{subfigure}{\linewidth}
		\centering
		\includegraphics[width=.5\linewidth]{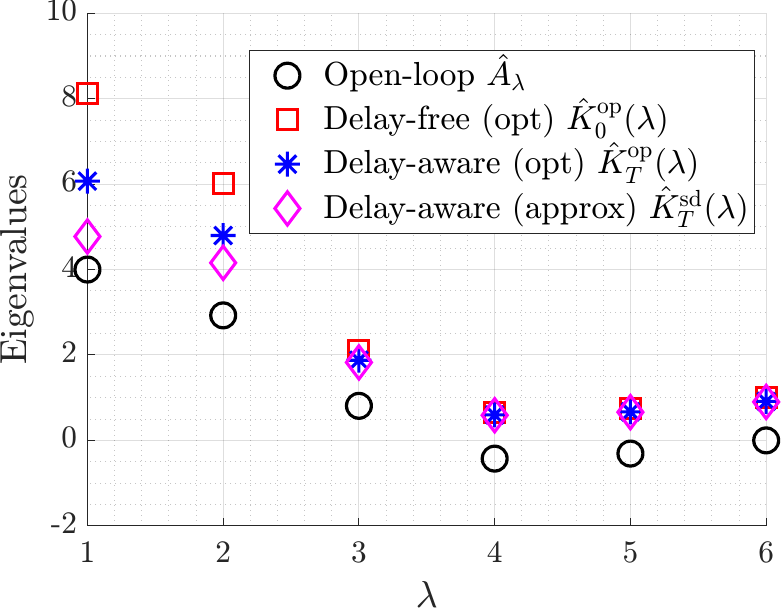}%
		\includegraphics[width=.5\linewidth]{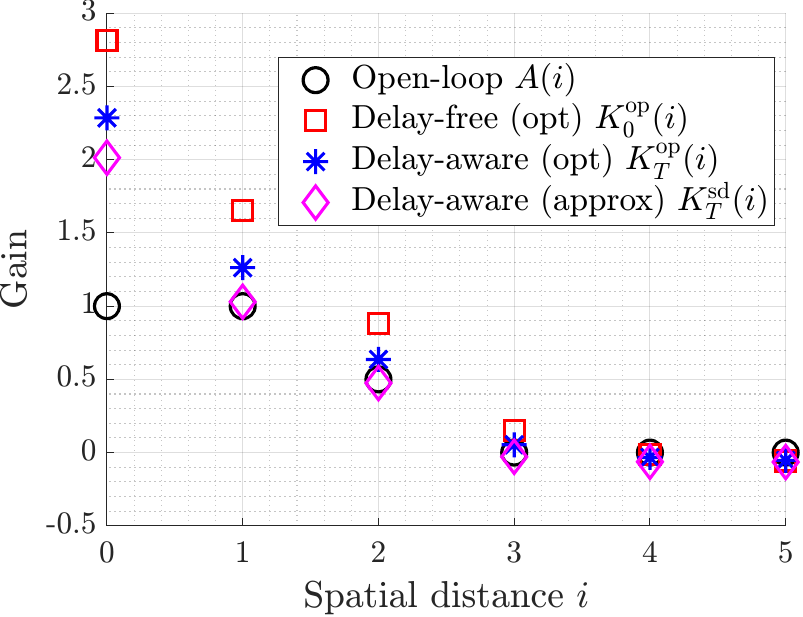}
		\caption{Delay $\delay=0.1$, control weight $r=1$.}
		\label{fig:Small_delays_r_1_T_0.1}
	\end{subfigure}\\
	\begin{subfigure}{\linewidth}
		\centering
		\includegraphics[width=.5\linewidth]{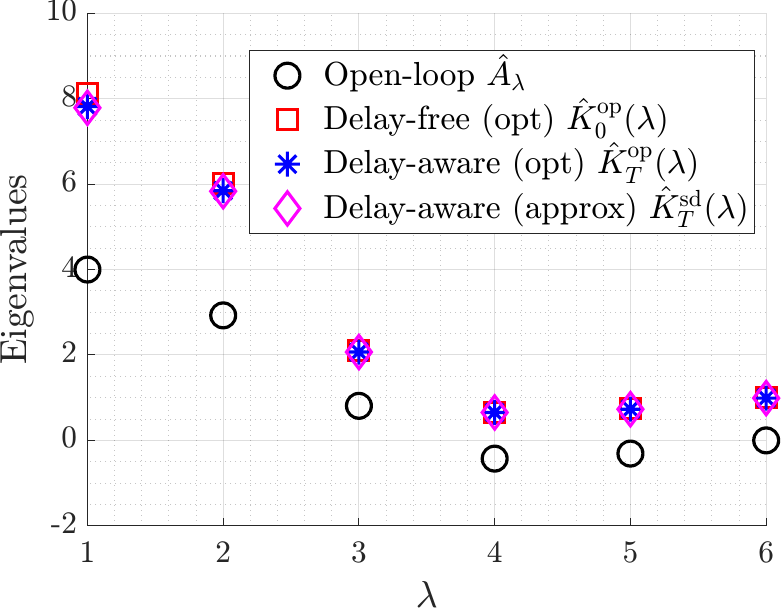}%
		\includegraphics[width=.5\linewidth]{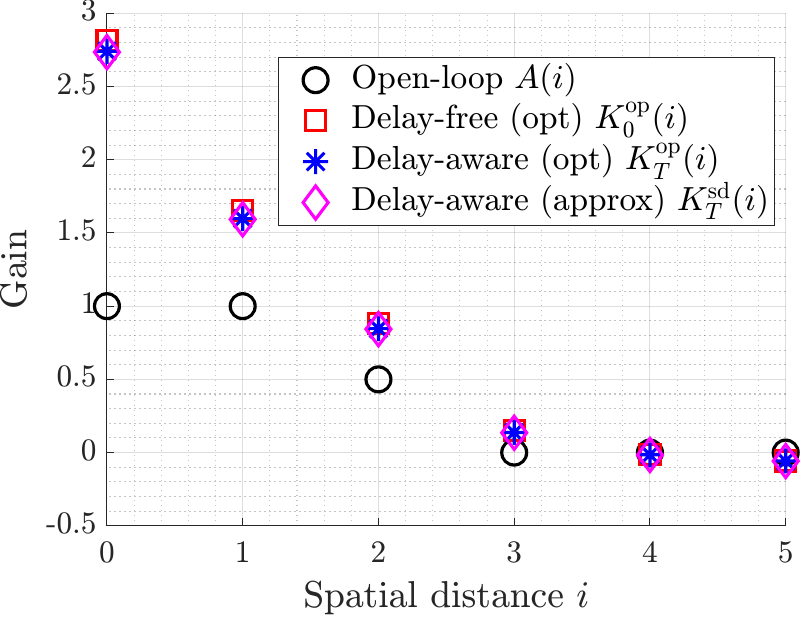}
		\caption{Delay $\delay=0.01$, control weight $r=1$.}
		\label{fig:Small_delays_r_1_T_0.01}
	\end{subfigure}\\
	\begin{subfigure}{\linewidth}
		\centering
		\includegraphics[width=.5\linewidth]{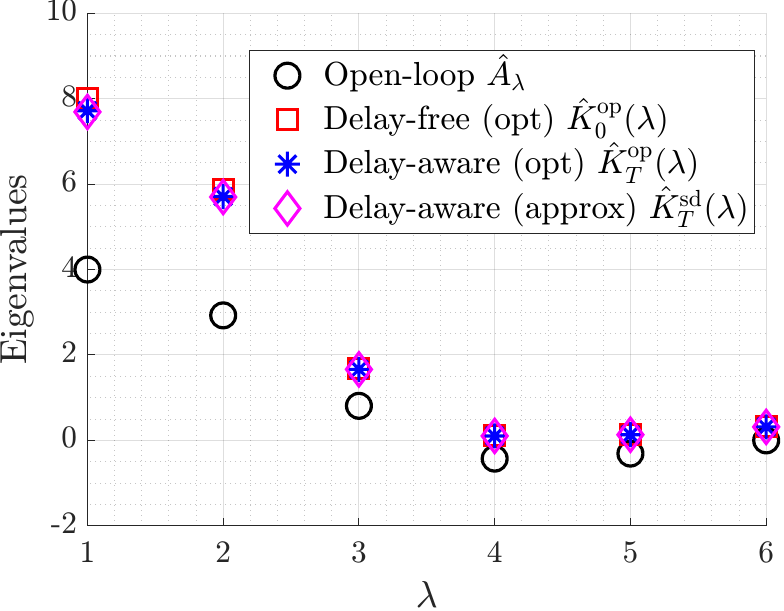}%
		\includegraphics[width=.5\linewidth]{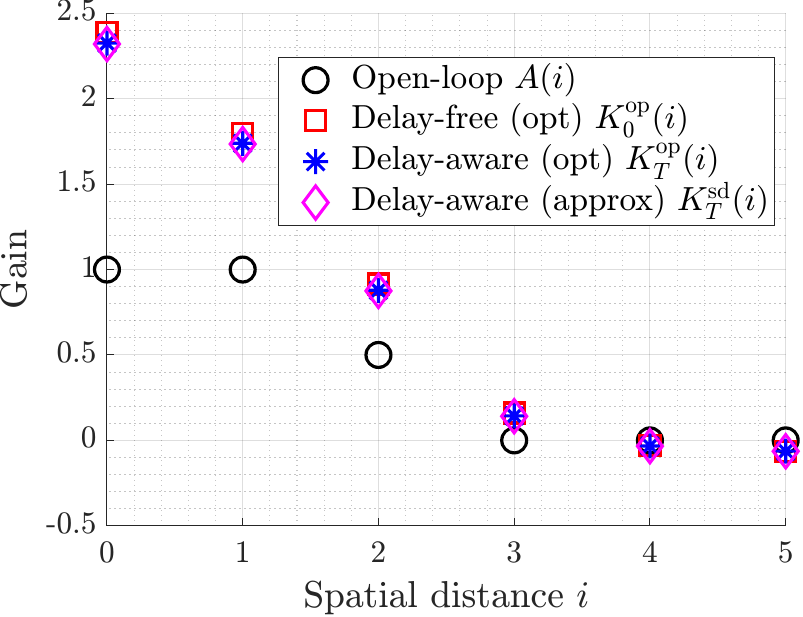}
		\caption{Delay $\delay=0.01$, control weight $r=10$.}
		\label{fig:Small_delays_r_10_T_0.01}
	\end{subfigure}
	\caption{Eigenvalues (left boxes) and gains (right boxes, both open loop and feedback) under small delays for system~\eqref{eq:circulant-dynamics} with $N=10$.}
	\label{fig:small-delays}
\end{figure}

\Cref{fig:small-delays} illustrates eigenvalues,
open-loop couplings, 
and feedback gains of open-loop and feedback matrices associated with system~\eqref{eq:circulant-dynamics} for $N=10$ agents and open-loop coupling gains among agents given by the vector
\begin{equation}
	A = \begin{bmatrix}
		1 & 1 & 0.5 & 0 & 0 & 0 & 0 & 0 & 0.5 & 1
	\end{bmatrix}.
\end{equation}
We represent only half-plus-one eigenvalues and gains in the plots because of symmetry.
The small-delay approximation~\eqref{eq:circulant-small-delay-gains} has a nonnegligible mismatch with the optimal feedback gains with delay $\delay=0.1$,
as shown in the right box of \autoref{fig:Small_delays_r_1_T_0.1},
while for $\delay=0.01$ they are almost identical to the numerically obtained optimal ones,
besides being similar to the delay-free optimal feedback gains.
All approximate and optimal feedback gains under delays are smaller than the corresponding optimal gains in the delay-free case,
and longer delays induce smaller gains.
Also,
increasing the control weight from $r=1$
to $r=10$ mostly affects the self-gains $\Kopk[\opt][\delay](0)$ and $\Kopk[\text{sd}][\delay](0)$,
which decrease from $2.8$ in \autoref{fig:Small_delays_r_1_T_0.01} to $2.3$ in \autoref{fig:Small_delays_r_10_T_0.01},
and the nearest-neighbor gains,
which increase from $1.6$ to $1.7$.

\section{Conclusion and Outlook}\label{sec:conclusion}

We study optimal feedback controllers for classes of spatially invariant systems where the measurements used for feedback control are subject to communication delays. These render the closed-loop dynamics a time-delay system.
First, we provide analytical and numerical insight into optimal proportional controllers for scalar delay systems.
Building upon these insights, we analyze the optimal control problem for spatially invariant systems with distributed and delayed measurements in two asymptotic cases: 
(i) expensive control regime $r\rightarrow+\infty$ and (ii) small delays $\delay\rightarrow0^+$,
characterizing spatial properties of the optimal feedback gain and its differences with the optimal gain for the standard delay-free scenario.
We also analyze two examples with a reaction-diffusion process over the real line and a multi-agent system in circulant interconnection. 
These corroborate our general theoretical result, 
illustrating how communication delays influence the optimal controllers' spatial spread.

This work makes a step forward towards understanding the structure of optimal controllers under communication delays,
especially regarding their spatial locality properties.
It suggests interesting follow-up research directions.
For example,
other delay models could be addressed,
or a more general formulation of the controller -- where the feedback operator need not be a time-independent spatial convolution.
The connection between communication delay and distributed controller architecture,
such as delays that depend on the communication topology~\cite{ballotta2023tcns,gupta2010delay},
should be further investigated to assess practical limitations of distributed control under real-world conditions.
	

	\appendix
	\numberwithin{equation}{subsection}

\subsection{Background on Spatially Invariant Systems}
\label{app:preliminaries}

Consider functions on space domain $\spacedom$ with spatial coordinate $x\in\spacedom$,
with $\mathcal{X}$ a locally compact abelian group.
The dual group of $\spacedom$ is $\Lambda$ and represents spatial frequency domain.

\begin{definition}[Translation invariance]
	Let $\mathcal{T}_x$ denote a translation operator defined as $[\mathcal{T}_x f](y) \doteq f (y - x)$ for any $f:\spacedom\to\spacedom$ and $x\in\spacedom$.
	An operator $\mathcal{A}$ is translation invariant if $\mathcal{T}_x \mathcal{A} = \mathcal{A}\mathcal{T}_x$ for every translation $\mathcal{T}_x :D(\mathcal{A}) \to D(\mathcal{A})$.
\end{definition}

\begin{definition}[Multiplication operator]
	A multiplication operator $\mathcal{M}$ is defined as $[\mathcal{M}f](x) \doteq M(x)f(x) \,\forall f\in\mathcal{D}(\mathcal{M})$,
	where $M$ is a measurable function called the \emph{symbol} of $\mathcal{M}$.
\end{definition}

\begin{definition}[Spatial Fourier transform]
	Let $f(\cdot,t)\in L^2(\spacedom)$ be a time-indexed function supported on $\spacedom$.
	We define the spatial Fourier transform pair of $f$ as follows:
	\begin{gather}
		\hat{f}_\lambda(t) = \fourier{f(\cdot,t)} \doteq \dfrac{1}{\sqrt{2\pi}}\int_{\spacedom} f(x,t) \e^{-jx\lambda}\de x\\
		f(x,t) = \ifourier[x]{\hat{f}_\lambda(t)} \doteq \dfrac{1}{\sqrt{2\pi}}\int_{\Lambda} \hat{f}_\lambda(t) \e^{jx\lambda}\de \lambda.
	\end{gather}
\end{definition}

The spatial Fourier transform \textit{diagonalizes} spatially invariant operator $\mathcal{A}$,
\ie the image $[\mathcal{A}f]$ of function $f$ through $\mathcal{A}$ transforms in spatial frequency domain to the multiplication $\hat{A}_\lambda \hat{f}_\lambda, \lambda\in\Lambda$.
The symbol $\hat{A}_\lambda$ of the multiplication operator in spatial frequency domain is called the \textit{Fourier symbol} of $\mathcal{A}$.

\subsection{Proof of \cref{prop:upper-bound}}\label{app:upper-bound}
{
	\renewcommandx{\Kopfc}[2][1={},2={}]{k_{#2}^{#1}}
	\renewcommand{\Aopf}{a}
	\renewcommand{\Aopftot}[1][]{a^{#1}}
    \renewcommand{\ub}{k_\text{u}}
    
	We use the implicit function theorem to prove that $\ub$ is decreasing with $\Aopf$ and convex in $\Aopf$.
	For $\Kopfc>|\Aopf|$,
	let
	\begin{equation}
		F(\Aopf,\Kopfc) \doteq \delay \sqrt{\Kopfc[2]-\Aopf^2} - \arccos\lr\Aopf\Kopfc[-1]\rr
	\end{equation}
    such that $F(\Aopf,\ub)=0$.
	The partial derivative $F_k(\Aopf,\Kopfc)$ is zero only if $\Kopfc[2]\delay = \Aopf$,
	which does not correspond to solutions of $F(\Aopf,\Kopfc) = 0$.
	This trivially holds true if $\Aopf\le0$.
	For $\Aopf>0$,
	assume $\ub^2\delay = \Aopf$,
	then $F(\Aopf,\ub) = 0$ reduces to
	\begin{equation}\label{eq:upper-bound-ineq-F}
		F(\ub^2\delay,\ub) = \ub\delay \sqrt{1 - \ub^2\delay^2} - \arccos\lr\ub\delay\rr \stackrel{\eqref{eq:upper-bound}}{=} 0
	\end{equation}
	for $\ub\delay \le1$.
	The derivative of the LHS of~\eqref{eq:upper-bound-ineq-F} w.r.t. $\ub$ is
	\begin{equation}
		F'(\ub^2\delay,\ub) = 2\delay\sqrt{1-\ub^2\delay^2} > 0 \quad \forall \; \ub\delay \in (0,1),
	\end{equation}
	hence $F(\ub^2\delay,\ub)$ increases with $\ub$ and the only real solution of~\eqref{eq:upper-bound-ineq-F} is $\ub=\nicefrac{1}{\delay}$,
    which is not valid because it implies $\ub=\Aopf$.
    Then,
    $\ub^2\delay \neq \Aopf$ and the implicit function theorem states that,
	for suitable intervals $\mathcal{I}_{\Aopf}$, $\mathcal{I}_{\Kopfc}$ and every $\Aopf\in\mathcal{I}_{\Aopf}$,
	the solution of $F(\Aopf,\Kopfc) = 0$ is given by $\ub=\ub(\Aopftot)\in\mathcal{I}_{\Kopfc}$.
	Moreover,
	for all $\Aopf$ such that $\Aopftot\ub^{-1}<\cos(\Aopftot\ub^{-1})$,
	it holds
	\begin{equation}\label{eq:upper-bound-ineq}
		\ub\delay > \delay \sqrt{\ub^2-\Aopf^2} \stackrel{\eqref{eq:upper-bound}}{=} \arccos\lr\Aopftot\ub^{-1}\rr > \Aopftot\ub^{-1}.
	\end{equation}
	Then,
	by continuity of $\ub(\Aopftot)$,
    it holds $\ub^{2}\delay > \Aopftot \; \forall \Aopf\in\mathcal{I}_\Aopf$.
	The derivative of $\ub(\Aopf)$ is
	\begin{equation}\label{eq:upper-bound-der}
		\ub'(\Aopf) = -\dfrac{F_{\Aopf}(\Aopf,\ub(\Aopf))}{F_{\Kopfc}(\Aopf,\ub(\Aopf))} = \dfrac{\ub(\Aopf)(\Aopf\delay - 1)}{\ub^2(\Aopf)\delay - \Aopf}.
	\end{equation}
	In view of $\ub>|\Aopf|\ge0$,
	the stability condition $\Aopf\delay - 1 < 0$,
	and $\ub^{2}\delay > \Aopftot$,
	it follows $\ub'(\Aopf) < 0$ for all $\Aopf\in\mathcal{I}_{\Aopf}$.
	As for the second derivative,
	standard calculations on~\eqref{eq:upper-bound-der} yield
	\begin{equation}
		\ub''(\Aopf) \propto -\ub'(\Aopf) \lr \delay + \dfrac{\Aopf}{\ub^2(\Aopf)}\rr > 0 \ \forall \Aopf\in\mathcal{I}_{\Aopf}
	\end{equation}
	where the term between brackets is always positive,
	by using the same argument of~\eqref{eq:upper-bound-ineq} for the case $\Aopf<0$.
	The rightmost limit in~\eqref{eq:upper-bound-limit} follows by continuity of $F$ and $F(\nicefrac{1}{\delay}, \nicefrac{1}{\delay}) = 0$~\cite{KuchlerLangevinEqs}.
	The leftmost limit in~\eqref{eq:upper-bound-limit} can be inferred by evaluating the derivative~\eqref{eq:upper-bound-der} at the limit.
	Let us assume that the limit is correct,
	which translates to the ansatz $\ub(\Aopf) \sim |\Aopf|$ for $\Aopf\rightarrow-\infty$.
	By continuity of~\eqref{eq:upper-bound-der},
	it follows
	\begin{equation}
		\lim_{\Aopf\rightarrow-\infty} \ub'(\Aopf) = \lim_{\Aopf\rightarrow-\infty} \dfrac{|\Aopf|(\Aopf\delay - 1)}{|\Aopf|(|\Aopf|\delay + 1)} = -1,
	\end{equation}
	consistently with the ansatz.
	We now rule out that at the limit $\ub(\Aopf) \to |\Aopf| + \kappa$ with $\kappa > 0$.
	Consider $F(\Aopf, |\Aopf| + \kappa)$.
	In view of $\delay>0$,
	it holds
	\begin{equation}
		\lim_{\Aopf\rightarrow-\infty} F(\Aopf, |\Aopf|+\kappa) = \lim_{\Aopf\rightarrow-\infty} \delay \sqrt{2|\Aopf|\kappa} \neq 0.
	\end{equation}
	Hence,
	it must be $\kappa = 0$ which validates the ansatz.
}

\subsection{Proof of \cref{prop:optimal-gain-tails}}\label{app:optimal-gain-small-a}

We first show that,
at the limit for $\Aopf\rightarrow-\infty$,
$f_\lambda$ tends to a constant function of $\Kopfc$ for $\Kopfc \in o(|\Aopf|)$ and it grows unbounded if $\Kopfc\in O(|\Aopf|)$.
At the limit for $\Aopf\rightarrow-\infty$,
the domain of $f_\lambda$ reduces to $(\Aopf, |\Aopf|]$ from~\cref{prop:upper-bound}.
Hence,
we restrict the analysis to the first two cases in~\eqref{eq:time-integral-function}.

\subsubsection{Case $\Kopfc\in(\Aopf, |\Aopf|)$}
We split this analysis based on the relationship between the asymptotic growth of $\Kopfc$ and $\Aopf$.

\paragraph{$\Kopfc\in O(|\Aopf|)$}
In view of $\Kopfc\in(\Aopf,|\Aopf|)$, 
we write $\Kopfc \sim \kappa|\Aopf|$ where $\kappa \in (0,1) $ is constant. 
It follows:
\begin{equation}\label{eq:f-approx-tail}
	f_\lambda(\kappa) \sim \frac{g(\kappa)}{2|\Aopftot|\sqrt{1-\kappa^2}}
\end{equation}
where we make dependence on $\kappa$ explicit and define
\begin{equation}
	g(\kappa)\doteq\dfrac{\kappa\sinh\lr|\Aopftot|\delay\sqrt{1-\kappa^2}\rr + \sqrt{1-\kappa^2}}%
	{1+\kappa\cosh\lr|\Aopftot|\delay\sqrt{1-\kappa^2}\rr}.
\end{equation}
We now split the function $g$ as follows:
\begin{equation}\label{eq:g-decomposition}
	g(\kappa) = g_1(\kappa) + g_2(\kappa)
\end{equation}
with
\begin{align}
	g_1(\kappa)	&\doteq\dfrac{\sqrt{1-\kappa^2}}%
	{1+\kappa\cosh\lr|\Aopftot|\delay\sqrt{1-\kappa^2}\rr}\label{eq:g1}\\
	g_2(\kappa)	&\doteq\dfrac{\kappa\sinh\lr|\Aopftot|\delay\sqrt{1-\kappa^2}\rr}{1+\kappa\cosh\lr|\Aopftot|\delay\sqrt{1-\kappa^2}\rr}.\label{eq:g2}
\end{align}
Next,
we approximate the hyperbolic functions.
To this aim,
we use the following limits that hold for $x\rightarrow\infty$:
\begin{align}
	\dfrac{1}{\sinh(x)} &\sim 2\e^{-x} \label{eq:appr-sinh}\\
	\dfrac{1}{1+\alpha\cosh(x)} &\sim \dfrac{2\e^{-x}}{\alpha+2\e^{-x}}\label{eq:appr-cosh}.
\end{align}
Combining~\eqref{eq:g1} with~\eqref{eq:appr-cosh} yields
\begin{equation}\label{eq:g1-approx}
	g_1(\kappa) \sim \dfrac{2\e^{-|\Aopftot|\delay\sqrt{1-\kappa^2}}\sqrt{1-\kappa^2}}{\kappa+2\e^{-|\Aopftot|\delay\sqrt{1-\kappa^2}}}.
\end{equation}
Combining~\eqref{eq:g2} with~\eqref{eq:appr-sinh}--\eqref{eq:appr-cosh} yields
\begin{equation}\label{eq:g2-approx}
	g_2(\Kopfc) \sim\dfrac{\kappa}{\kappa + 2\e^{-|\Aopftot|\delay\sqrt{1-\kappa^2}}}.
\end{equation}
Finally,
combining~\eqref{eq:g-decomposition} with~\eqref{eq:g1-approx}--\eqref{eq:g2-approx} yields
\begin{equation}
	g(\kappa) \sim \dfrac{\kappa + 2\e^{-|\Aopftot|\delay\sqrt{1-\kappa^2}}\sqrt{1-\kappa^2}}{\kappa + 2\e^{-|\Aopftot|\delay\sqrt{1-\kappa^2}}} \sim 1.
\end{equation}
From~\eqref{eq:f-approx-tail},
we obtain $f_\lambda(\Kopfc) \sim \frac{1}{2|\Aopf|\sqrt{1-\kappa^2}}$ and the cost function evaluates
$J_\lambda(\Kopfc) \sim \frac{1 + r\kappa^2|\Aopf|^2}{2|\Aopf|\sqrt{1-\kappa^2}}$.
This value tends to infinity as $\Aopf\rightarrow-\infty$ for fixed values of $r$ and $\kappa$.

\paragraph{$\Kopfc\in o(|\Aopf|)$}
We evaluate the limit of function $f_\lambda$ in analogous manner as done above.
It holds:
\begin{equation}\label{eq:f-approx-tail-1}
	f_\lambda(\Kopfc) \sim \frac{h(\Kopfc)}{2|\Aopftot|}, \quad 
	h(\Kopfc)\doteq\dfrac{\dfrac{\Kopfc}{|\Aopftot|}\sinh\lr|\Aopftot|\delay\rr + 1}%
	{1+\dfrac{\Kopfc}{|\Aopftot|}\cosh\lr|\Aopftot|\delay\rr}
\end{equation}
where $\ell\sim|\Aopftot|$ for $\Aopf\rightarrow-\infty$.
Similarly to case \textit{a)},
we get
\begin{equation}\label{eq:h-asympt}
	h(k) \sim \dfrac{1}{\dfrac{\Kopfc}{|\Aopftot|}\dfrac{1}{2\e^{-|\Aopftot|\delay}}+1} + \dfrac{1}{\dfrac{|\Aopftot|}{\Kopfc}2\e^{-|\Aopftot|\delay} + 1} \equiv 1.
\end{equation}
From~\eqref{eq:f-approx-tail-1}--\eqref{eq:h-asympt},
we get $f_\lambda(\Kopfc) \sim \frac{1}{2|\Aopf|}$ and the cost has limit $J_\lambda(\Kopfc) \sim \frac{1+r\Kopfc[2]}{2|\Aopf|}$.
This is minimized by $\Kopfc[\opt][\delay] = 0 \in o(|\Aopf|)$,
with corresponding limit of the cost $J_\lambda(\Kopfc[\opt][\delay])\rightarrow0$.

We conclude that the minimizer of $J_\lambda$ when $f_\lambda$ is defined by the first case in~\eqref{eq:time-integral-function}
has limit $\Kopfc[\opt][\delay] \rightarrow 0$ as $\Aopf\rightarrow-\infty$.

\subsubsection{Case $\Kopfc = |\Aopf|$}
As $\Aopf\rightarrow-\infty$, 
the second case yields $f_\lambda(|\Aopf|) \rightarrow \nicefrac{\delay}{4}$ and $J_\lambda(|\Aopf|) \sim (1 + r|\Aopf|^2)\nicefrac{\delay}{4} = +\infty$.

\subsubsection{All Cases $\Kopfc\in(\Aopf,|\Aopf|]$}
Noting that the limit of the cost $J_\lambda$ in the regime $\Kopfc=|\Aopf|$ is infinite for any $r>0$,
it follows that,
at the limit for $\Aopf\rightarrow-\infty$,
the minimum of $J_\lambda$ is achieved in the first case $\Kopfc\in(\Aopf,|\Aopf|)$ and the minimum point $\Kopfc[\opt][\delay]$ tends to zero for any positive $r$.

By virtue of the observation $\Kopfc[\opt][\delay]\rightarrow0$ and because the cost function is continuous and positive,
its argmin is also continuous~\cite{rockafellar1998variational} and we can compute the limit of the minimizer $\Kopfc[\opt][\delay]$ via a McLaurin expansion of $J_\lambda$.
Consider the quadratic approximation $Q_J$ such that $J_\lambda(\Kopfc) = Q_{J}(\Kopfc) + o(\Kopfc[2])$,
\begin{equation}\label{eq:obj-function-tail-taylor}
	\begin{aligned}
		Q_J(\Kopfc)	&\doteq J_\lambda(0) + J_\lambda'(0)\Kopfc + J_\lambda''(0)\dfrac{\Kopfc[2]}{2}\\
		&\propto 1 + \dfrac{\e^{\delay\Aopftot}}{\Aopftot}\Kopfc + \lr\dfrac{\e^{2\delay\Aopftot}+1}{2\Aopftot[2]} + r \rr\Kopfc[2].
	\end{aligned}
\end{equation}
Setting $Q_J'(\Kopfc)=0$,
for $r>0$ we obtain the following solution that tends to the minimizer of $J_\lambda(\Kopfc)$ as $\Aopf\rightarrow-\infty$:
\begin{equation}\label{eq:kstar-gen-tail}
	\Kopfc[\opt,Q_J][\delay] = \dfrac{1}{2r|\Aopftot|} \e^{\delay\Aopftot}.
\end{equation}

\subsection{Proof of~\cref{thm:optimal-gain-expensive-regime}}\label{app:optimal-gain-expensive-regime}

Let $\epsilon\doteq\nicefrac{1}{r}$.
We rewrite the cost function of~\eqref{eq:prob-optimal-control-decoupled-explicit} as
\begin{equation}\label{eq:obj-function}
	J_\lambda(\Kopfc) = r\bar{J}_\lambda(\Kopfc), \quad \bar{J}_\lambda(\Kopfc)\doteq\lr\epsilon + \Kopfc[2]\rr f_\lambda(\Kopfc).
\end{equation}
It holds $\arg\min_{k}J_\lambda(k) = \arg\min_{k}\bar{J}_\lambda(k)$.
Under the assumption $\Aopftot < 0$,
because $f_\lambda$ is positive and continuous,
the minimizer of $\bar{J}_\lambda$ tends to zero as $\epsilon\rightarrow0^+$~\cite{rockafellar1998variational}.
Hence,
analogously to~\eqref{eq:obj-function-tail-taylor},
we consider the quadratic approximation $Q_{\bar{J}_\lambda}$ given by the second-order McLaurin expansion of $\bar{J}_\lambda$:
\begin{equation}\label{eq:obj-function-taylor}
	Q_{\bar{J}_\lambda}(\Kopfc) \doteq \bar{J}_\lambda(0) + \bar{J}_\lambda'(0)\Kopfc + \bar{J}_\lambda''(0)\dfrac{\Kopfc[2]}{2}.
\end{equation}
By continuity of $\bar{J}_\lambda$,
its minimum point tends to the minimum point of $Q_{\bar{J}_\lambda}$
as $\epsilon\rightarrow0^+$~\cite{rockafellar1998variational}.
Setting the derivative of $Q_{\bar{J}_\lambda}$ equal to zero,
we obtain the following solution:
\begin{equation}\label{eq:kstar-expensive-approx-full}
	\Kopfc[\opt,Q_{\bar{J}_\lambda}][\delay](\lambda)
	=\dfrac{-\epsilon\Aopftot\e^{\delay\Aopftot}}{\epsilon\e^{2\delay\Aopftot} + \epsilon + 2\Aopftot[2]}.
\end{equation}
With $\Aopftot$ bounded away from zero,
we simplify the expression above via its Taylor expansion about $\epsilon=0$:
\begin{equation}\label{eq:kstar-expensive-approx-proof}
	\Kopfc[\opt,Q_{\bar{J}_\lambda}][\delay](\lambda) = \dfrac{\epsilon}{2|\Aopftot|}\e^{\delay\Aopftot} + o(\epsilon).
\end{equation}

\subsection{Proof of~\cref{prop:optimal-gain-small-delay}}\label{app:optimal-gain-small-delay}

The first-order Taylor expansion of the cost function in~\eqref{eq:prob-optimal-control-decoupled-explicit} about $\delay=0$ is
\begin{equation}\label{eq:variance-approx-small-T-general}
	J_\lambda(\Kopfc,\delay) 
	=\dfrac{1+\Kopfc\delay}{2(\Kopfc-\Aopftot)}(1+r\Kopfc[2]) + \tilde{J}_\lambda(\Kopfc,\delay)
\end{equation}
where $\tilde{J}_\lambda(\Kopfc,\delay) \in o(\delay)$ and is a rational function of $\Kopfc$.
The Moore-Osgood theorem implies that $\frac{\de\tilde{J}_\lambda(\Kopfc,\delay)}{\de\Kopfc}\in o(\delay)$.
The roots of a polynomial vary continuously with the coefficients of the polynomial,
including vanishing leading coefficients,
whereby new roots emerge from a neighborhood of the point at infinity~\cite{rootsContinuity}.
Hence,
in view of the bounded domain $\Kopfc\in(\Aopf,\ub)$,
the critical points of the first-order approximation of $J_\lambda$ w.r.t. $\delay$ in~\eqref{eq:variance-approx-small-T-general} continuously approximate the relevant (bounded) critical points of $J_\lambda$ as $\delay\to0$.
Setting the derivative of $J_\lambda$ w.r.t. $\Kopfc$ equal to zero and neglecting $\tilde{J}_\lambda$ yields
\begin{equation}\label{eq:kstar-small-T-general-eqn}
	2r\delay \Kopfc[3] + (1-3\delay\Aopftot)r\Kopfc[2] - 2r\Aopftot \Kopfc - 1 - \delay\Aopftot = 0.
\end{equation}
Equation~\eqref{eq:kstar-small-T-general-eqn} is cubic and can be solved in closed form through Cardano's formula.
The three complex solutions are
\begin{equation}\label{eq:cubic-roots}
	\Kopfc[\text{c}] = - \dfrac{b}{3a} - \dfrac{1}{3a} \lr C + \dfrac{\Delta_0}{C} \rr
\end{equation}
with
\begin{align}
	a &= 2r\delay  \qquad 
	b = (1-3\delay\Aopftot)r \label{eq:ab} \\
	c &= - 2r\Aopftot \qquad
	d = - 1 - \delay\Aopftot\\
	\Delta_0 	&= b^2 - 3ac= r^2(1 + 3\delay\Aopf)^2\label{eq:delta0}\\
	\Delta_1 	&\begin{aligned}[t]
					&= 2b^3 - 9abc + 27a^2d\\
							&= \begin{multlined}[t]
								r^3(2 + 18\delay\Aopf - 54\delay^2\Aopf^2 - 54\delay^3\Aopf^3)\\
								-108r^2\delay^2(\delay\Aopf + 1).
							\end{multlined}
				\end{aligned}\label{eq:delta1}\\
	C &= \sqrt[3]{\dfrac{\Delta_1 + \sqrt{\Delta_1^2 - 4\Delta_0^3}}{2}}\label{eq:C}
\end{align}
The cubic root in~\eqref{eq:C} represents the three complex roots of the argument and
only one of these is real and positive by Descartes' rule of signs,
hence corresponding to the unique minimum point of $J_\lambda$.
Let $C^3 = r^3\rho\e^{j\theta}$.
We approximate $\rho$ and $\theta$ to explicitly compute the term between brackets in~\eqref{eq:cubic-roots} that corresponds to the unique real positive solution of~\eqref{eq:kstar-small-T-general-eqn}.
From~\eqref{eq:C} and $\Delta_1^2 - 4\Delta_0^3 < 0$,
we compute $\rho$ as
\begin{equation}\label{eq:rho}
	\rho	= \dfrac{1}{r^3}\sqrt{\lr\dfrac{\Delta_1}{2}\rr^{\!\!2} + \dfrac{-\lr\Delta_1^2 - 4\Delta_0^3\rr}{4}} = \sqrt{\dfrac{\Delta_0^3}{r^6}}
			= \lr1+3\delay\Aopf\rr^{\!3}
\end{equation}
and,
by $\delay\Aopf\rightarrow0$,
we approximate $\theta$ via Taylor expansion as
\begin{equation}\label{eq:theta}
	\theta	= \arctan\dfrac{\sqrt{4\Delta_0^3 - \Delta_1^2}}{\Delta_1} =\delay\sqrt{108\gamma}\lr1-4\delay\Aopf\rr + o(\delay^2)
\end{equation}
where we define
\begin{equation}
	\gamma \doteq \Aopf^2 + r^{-1}.
\end{equation}
We now evaluate the term between brackets in~\eqref{eq:cubic-roots} as
\begin{equation}\label{eq:cubic-roots-part}
	C + \dfrac{\Delta_0}{C} = r\lr\rho^{\frac{1}{3}} + \rho^{-\frac{1}{3}}\dfrac{\Delta_0}{r^2}\rr\cos\lr\dfrac{\theta+2\pi}{3}\rr.
\end{equation}
From~\eqref{eq:delta0} and~\eqref{eq:rho},
the first factor in~\eqref{eq:cubic-roots-part} amounts to
\begin{equation}\label{eq:rhos}
	\rho^{\frac{1}{3}} + \rho^{-\frac{1}{3}}\dfrac{\Delta_0}{r^2} = 2 + 6\delay\Aopf.
\end{equation}
From~\eqref{eq:theta},
we approximate the second factor in~\eqref{eq:cubic-roots-part} via Taylor expansion as
\begin{equation}\label{eq:cos}
	\begin{aligned}
		\cos\lr\dfrac{\theta+2\pi}{3}\rr 	&= \cos\dfrac{\theta}{3}\cos\dfrac{2\pi}{3} - \sin\dfrac{\theta}{3}\sin\dfrac{2\pi}{3}\\
											&=-\dfrac{1}{2}\lr1 - \dfrac{\theta^2}{18}\rr -\dfrac{\theta}{3}\dfrac{\sqrt{3}}{2} + o(\theta^2)\\
											&=\dfrac{1}{2} - 3\delay\sqrt{\gamma} + \lr3\gamma + 12\Aopf\sqrt{\gamma}\rr\delay^2 
												 + o(\delay^2).
	\end{aligned}
\end{equation}
Putting together~\eqref{eq:rhos}--\eqref{eq:cos},
we evaluate~\eqref{eq:cubic-roots-part} as
\begin{multline}\label{eq:C-2}
	C + \dfrac{\Delta_0}{C} = r \lr - 1 - 3\delay\Aopf - 6\delay\sqrt{\gamma} + 6\delay^2\Aopf^2 \right.\\
	\left.+ 6\delay^2\Aopf\sqrt{\gamma} + 6\dfrac{\delay^2}{r}\rr + o(\delay^2).
\end{multline}
Finally,
we plug~\eqref{eq:ab} and~\eqref{eq:C-2}
into~\eqref{eq:cubic-roots} to compute the real positive solution of~\eqref{eq:kstar-small-T-general-eqn} as
\begin{equation}\label{eq:kstar-small-T-proof}
	\Kopfc[\text{c},+] = \Aopftot + \sqrt{\gamma} - \lr\Aopftot\lr\Aopf+\sqrt{\gamma}\rr + \dfrac{1}{r}\rr\delay + o(\delay).
\end{equation}
We conclude that the optimal solution $\Kopfc[\opt][\delay](\lambda)$ that minimizes $J_\lambda$ in~\eqref{eq:variance-approx-small-T-general} is asymptotic to $\Kopfc[\text{c},+]$ as $\delay\rightarrow0$.

\subsection{Optimal Cost in Expensive Control Regime}\label{app:optimal-cost-expensive-regime}

Due to the integral form of the cost $J$,
we focus on the integrand $J_\lambda$.
For convenience of notation,
we use again $\epsilon\doteq\nicefrac{1}{r}$.
We first evaluate the function $f_\lambda(\Kopfc[\text{ex}][\delay](\lambda))$:
\begingroup
\allowdisplaybreaks
\begin{align}\label{eq:fstar-expensive}
	\begin{split}
		f_\lambda(\Kopfc[\text{ex}][\delay](\lambda)) 	&= \dfrac{-\Kopfc[\text{ex}][\delay](\lambda)\sinh\lr\ell\delay\rr - \ell}{2\ell\lr \Aopftot-\Kopfc[\text{ex}][\delay](\lambda)\cosh\lr\ell\delay\rr\rr}\\
		&\stackrel{(i)}{=} \dfrac{1}{2|\Aopftot|} \dfrac{4\Aopftot[2] + \epsilon\lr1-\e^{2\Aopftot\delay}\rr}{4\Aopftot[2] + \epsilon\lr1+\e^{2\Aopftot\delay}\rr}\\
		&\stackrel{(ii)}{=} \dfrac{1}{2|\Aopftot|} \lr1-\dfrac{\epsilon}{2\Aopftot[2]}\e^{2\Aopftot\delay}\rr + o(\epsilon)
	\end{split}
\end{align}
\endgroup
where $(i)$ follows by plugging in $\Kopfc[\text{ex}][\delay](\lambda)$ and $\ell\sim|\Aopftot|$
and $(ii)$ is a first-order Taylor approximation about $\epsilon=0$.
Then,
the optimal cost evaluates as follows:
\begin{equation}\label{eq:jstar-expensive-proof}
	\begin{aligned}
		J_\lambda\lr\Kopfc[\text{ex}][\delay](\lambda)\rr	&= r\lr\epsilon + \lr\Kopfc[\text{ex}][\delay](\lambda)\rr^2\rr f_\lambda\lr\Kopfc[\text{ex}][\delay](\lambda)\rr\\
								&\stackrel{(i)}{\sim} \dfrac{1}{2|\Aopftot|}\lr1+\dfrac{\epsilon\e^{2\Aopftot\delay}}{4\Aopftot[2]}\rr\lr1-\dfrac{\epsilon\e^{2\Aopftot\delay}}{2\Aopftot[2]}\rr\\
								&\stackrel{(ii)}{\sim} \dfrac{1}{2|\Aopftot|} \lr1-\dfrac{\epsilon}{4\Aopftot[2]}\e^{2\Aopftot\delay}\rr
	\end{aligned}
\end{equation}
where $(i)$ follows from plugging in $\Kopfc[\text{ex}][\delay](\lambda)$ and~\eqref{eq:fstar-expensive} without the term in $o(\epsilon)$
and $(ii)$ from neglecting further terms in $o(\epsilon)$.
	
	\if0\mode
	\input{Bio/bios}	
	\fi
	
\end{document}